\documentclass[a4paper, 11pt]{amsart}

\usepackage[english]{babel}
\usepackage{amsthm, amssymb}
\usepackage{color}
\usepackage{float}
\usepackage{algpseudocode}
\usepackage{xcolor}
\usepackage{amsmath}
\usepackage{url}
\usepackage{bm}
\usepackage[a4paper]{geometry}
\usepackage{booktabs}

\usepackage[inline,shortlabels]{enumitem}
\usepackage[pagebackref]{hyperref}
\hypersetup{
 colorlinks=true,
 linkcolor=blue,
 filecolor=blue, 
 urlcolor=blue,
 citecolor=blue,
}

\numberwithin{equation}{section}
\newtheorem{theorem}{Theorem}[section] 
\newtheorem{lemma}[theorem]{Lemma} 
\newtheorem{proposition}[theorem]{Proposition}
\newtheorem{prop}[theorem]{Proposition}
\newtheorem{corollary}[theorem]{Corollary}
\newtheorem{definition}[theorem]{Definition}
\newtheorem{defn}[theorem]{Definition}
\theoremstyle{definition}
\newtheorem{remark}[theorem]{Remark}

\DeclareMathOperator{\GL}{GL}
\DeclareMathOperator{\SL}{SL}
\DeclareMathOperator{\SU}{SU}

\DeclareMathOperator{\im}{im}

\DeclareMathOperator{\dist}{dist}
\DeclareMathOperator{\diag}{diag}
\DeclareMathOperator{\Lip}{Lip}
\DeclareMathOperator{\Span}{span}

\newcommand{\Fa}{F_a}
\newcommand{\la}{\lambda}
\newcommand{\ra}{\rightarrow}

\newcommand{\R}{\mathbb{R}}

\renewcommand{\phi}{\varphi}
\renewcommand{\leq}{\leqslant}
\renewcommand{\geq}{\geqslant}
\renewcommand{\epsilon}{\varepsilon}

\newcommand{\ol}{\overline}

\newcommand{\x}{x^{-1}}
\newcommand{\y}{y^{-1}}
\newcommand{\xx}{x^{-1}} 
\newcommand{\yy}{y^{-1}}

\newcommand{\Tr}{\mathrm {Trace}}

\newcommand{\N}{\mathbb N}

\newcommand{\Z}{\mathbb Z}
\newcommand{\Q}{\mathbb Q}
\newcommand{\C}{\mathbb C}
\renewcommand{\P}{\mathbb P}

\makeatletter
\@namedef{subjclassname@2020}{%
	\textup{2020} Mathematics Subject Classification}
\makeatother

\begin{document}

\title{Constructing free subgroups of matrix groups}

\author{Emmanuel Breuillard} 
\address[Emmanuel Breuillard]{
Mathematical Institute,
University of Oxford, 
Oxford OX2 6GG,
United Kingdom}

\author{Bettina Eick}
\address[Bettina Eick]{
Institute for Analysis and Algebra,
Technische Universit\"at Braunschweig,
38106 Braunschweig, Germany}

\author{Tommy Hofmann}
\address[Tommy Hofmann]{
University of Siegen, 57068 Siegen, Germany}

\author{Derek Holt}
\address[Derek Holt]{
Mathematics Institute,
University of Warwick,
Coventry CV4 7AL,
United Kingdom}

\author[E.\ A.\ O'Brien]{E.\ A.\ O'Brien}
\address[E.\ A.\ O'Brien]{
Department of Mathematics,  
University of Auckland,
Private Bag 92019,
Auckland, 
New Zealand }

\renewcommand{\shortauthors}{E.\ Breuillard, B.\ Eick, T.\ Hofmann, 
D.\ Holt, E.A.\ O'Brien}


\thanks{
O'Brien was supported by the Marsden Fund of New Zealand Grant 23-UOA-080
and by a Research Award of the Alexander von Humboldt Foundation.
Hofmann gratefully acknowledges support by the Deutsche
Forschungsgemeinschaft --- Project-ID 286237555 - TRR 195; and
Project-ID 539387714.  We thank the Institute for Analysis and Algebra 
at the Technische Universit\"at Braunschweig for its hospitality and 
support of associated research visits.}

\begin{abstract}
The Tits alternative asserts that a finitely generated linear group 
is either solvable-by-finite or contains a non-abelian free subgroup.
Given a non-solvable-by-finite subgroup $G$ of 
$\GL(d,F)$, where $F$ is the field of complex numbers or a non-archimedean 
local field, we present methods to construct explicitly a non-abelian 
2-generator free subgroup of $G$ 
in specific cases. Our implementations in 
the computer algebra system Oscar are publicly available. 
\end{abstract}

\maketitle

\section{Introduction}

The Tits alternative \cite{tits} asserts that a finitely generated linear
group $G$ is either solvable-by-finite or contains a non-abelian free 
subgroup. This theorem partitions linear groups into two classes which 
require separate treatment. Algorithms of Beals \cite{Beals99,Beals01} 
and Detinko, Flannery \& O'Brien \cite{DFO11} allow us to decide 
if $G$ is solvable-by-finite. 

We introduce the first explicitly quantitative version of the 
Tits alternative and, 
based on this, the first practical algorithms to construct an explicit 
non-abelian free subgroup of a finitely generated linear group $G$ 
defined over the complex field or over a non-archimedean local field.  
To do so, we extend the techniques of Tits 
\cite{tits}, and build on  
later results of Dixon \cite{Dix88}, Breuillard \cite{Bre08}, 
and Breuillard and Gelander \cite{breuillard-gelander, BG-invent}.

Our methods use randomization. Each method takes as input a generating set 
for $G$. It may or may not terminate. If it terminates, then it returns a 
{\em free pair} for $G$: this is a free generating set of a non-abelian free 
subgroup of $G$. The elements in the free pair are given as words in the 
generators of $G$. We exhibit conditions on $G$ which 
guarantee termination.

\medskip
\noindent 
{\bf Linear groups over the field of complex numbers.}
We introduce three methods to construct a free pair for $G \leq \GL(d,\C)$: 
the {\it bi-proximal method} (Theorem \ref{Ebound}); its 
variant called the {\it conjugator method} (Theorem \ref{Bbound-alternative}); 
and the {\it singular value decomposition method} 
(Theorem \ref{SVDtheorem0}).  
These methods rely on the existence of certain types of 
elements: $g \in G$ is \emph{proximal} if its eigenvalues 
$\la_1(g), \ldots, \la_d(g)$ sorted by absolute values 
satisfy $|\la_1(g)| > |\la_2(g)|$, and $g$ 
is \emph{bi-proximal} if both $g$ and $g^{-1}$ are \emph{proximal}. We show in 
Theorem \ref{termC}(1) that the methods terminate if $G$ acts 
strongly irreducibly on $\C^d$ and contains a proximal element. 
These conditions are satisfied if, for example, $G$ is Zariski-dense in 
$\SL(d, \R)$, see Theorem \ref{termC}(2).

\medskip
\noindent 
{\bf Linear groups over non-archimedean local fields.}
We introduce variations of the conjugator method and the singular value
decomposition method for 
groups $G \leq \GL(d,K)$, where $K$ is a non-archimedean local field; see 
Theorems \ref{Bboundultra} and \ref{SVDtheorem}.
These methods rely on the existence
of {\em generalized} bi-proximal elements. Their termination is discussed
in Theorem \ref{termK}; Zariski-density in $\SL(d, K)$ and
unboundedness are required.

\medskip
\noindent 
{\bf Linear groups over the field of rational numbers.}
Of special interest is the case where $G \leq \GL(d, \Q)$, since $\Q$ 
embeds in both $\C$ and $\Q_p$,
the field of $p$-adic numbers, for all primes $p$. 
Certain conditions must be satisfied for termination. 
For example, the
conjugator method over $\Q_p$ requires that $G$ is not conjugate to a
subgroup of $\GL(d,\Z_p)$. 
This condition can hold only for the finitely many primes $p$ dividing a
denominator of an entry of a generator or its inverse;
Theorem~\ref{bounded} allows us to decide whether it holds for a given $p$.

\medskip
\noindent
A major challenge facing the implementations over $\C$ is to ensure
that the necessary delicate arithmetic is carried out accurately, avoiding
accumulation of rounding errors. We chose the computational
algebra system Oscar \cite{OSCAR-book, Hecke} as our implementation
platform, since it supports {\it interval arithmetic}
\cite{AH1983, Johansson2013arb}, so addressing this concern.

The outline of the paper is the following. In Section \ref{sec_prelim},
we review some background, in particular the Ping-Pong Lemma of Tits, and 
introduce notation that is used throughout. 
In Sections \ref{sec_biprox}--\ref{sec-padicSVD}, we develop
the theory necessary for our methods.
In Section \ref{algo_var1}, we describe 
the resulting algorithms
and report on their publicly available implementations in Oscar. 
In Section \ref{section_cones}, we describe a procedure
which uses cones to verify that two rational matrices 
form a free pair.
In Section \ref{sec_burau}, we report on an application 
to the 4-strand Burau representation. 
Finally, we comment briefly on the performance of our Oscar implementations
which are available as a GitHub package \cite{GITPAGE}.

We note that we do not resolve the problem of 
deciding whether or not two given invertible matrices generate a 
non-abelian free subgroup. 
It seems likely that this problem 
is \emph{undecidable}; see \cite{breuillard-kocharyan} where 
the related non-commutative Skolem problem is shown to be undecidable for 
subgroups of $\GL(d, \Q)$ when $d \ge 9$.

\medskip
\noindent \emph{AI use.} This paper was developed and written
by humans. LLMs (ChatGPT and Claude) were used 
to provide a sharp multiplicative constant in Lemma~\ref{conjlem},
for proofreading, and to produce Oscar scripts.

\section{Preliminaries and notation}
\label{sec_prelim}

Let $F$ denote $\C$ or a non-archimedean local field $K$.
For $\alpha \in F$, denote its absolute value by $|\alpha|$. 
If $v \in F^d$ and $A \in M(d,F)$, then $vA$ denotes 
multiplication of the row vector $v$ by the matrix $A$.
Let $\la_1(A), \ldots, \la_d(A)$ denote the roots of the characteristic 
polynomial of $A$ sorted by their absolute values, where 
\[ \Lambda(A) := \max \{ |\la_i(A)| \, : \, 1 \leq i \leq d \}.\]

\subsection{Norms and distance over $\C$}
\label{sec_biprox:distance}
We establish notation for well-known concepts.

Let $u = (u_1, \ldots, u_d)$ and $v = (v_1, \ldots, v_d)$ be vectors 
in $\C^d$. 
\begin{itemize}
\item
$||u|| = (\sum_{i=1}^d |u_i|^2)^{1/2}$.
\item
$u \wedge v$ is the vector with entries $u_i v_j - u_j v_i$ for $i < j$.
\item $u \cdot v = \sum_{i=1}^d u_i \overline{v}_i$.
\item
If $u,v \neq 0$, then $\dist(u,v) := || u \wedge v || / (||u|| \cdot ||v|| )$.
\end{itemize}

If $V \leq \C^d$ is a hyperplane with normal vector $v$ and 
$u \not=0$, then  
\[ \dist(u,V) := \inf_{w \in V\setminus\{0\}} \dist(u,w) 
              = |u \cdot v | / (||u|| \cdot ||v||).\]

Let $A = (A_{i,j}) \in M(d,\C)$. 
Recall that its adjoint $A^*$  is the conjugate transpose $\ol{A^T}$. 
The operator norm of $A$ is $\|A\|_o = 
\max_{\|v\| = 1} \|vA\| 
= \sqrt{\Lambda(A^* A)}$.
Observe that $\Lambda(A) \leq ||A||_o$.
In Section \ref{sec_untra}, we introduce analogous concepts 
for a non-archimedean local field.

\subsection{The Ping-Pong Lemma}
We first recall (an abstract version of) 
the lemma as presented by Tits \cite[Proposition 1.1]{tits}. 

\begin{lemma}\label{Tits-version}
Suppose a group $G$ acts on a set $P$. Let $p\in P$ be a point, 
let $X$ and $Y$ be subsets of $P$, and let 
$\langle x \rangle$ and $\langle y \rangle$ be infinite cyclic subgroups of 
$G$. If
\begin{itemize}
\item $p \notin X \cup Y$,
\item $(X \cup \{p\})g \subset Y$ for all 
$g \in \langle y \rangle \setminus\{1\}$,
       and
\item  $(Y \cup \{p\})g \subset X$ for all $g \in \langle x \rangle \setminus\{1\}$,
\end{itemize}
then $x$ and $y$ form a free pair.
\end{lemma}

As remarked by Tits, conversely, if $x$ and $y$ form a free pair, then let $\langle x,y \rangle$ act on itself and let $X$ (resp.\ $Y$) be the set of non-trivial words ending with the letter $x$ or $x^{-1}$ (resp.\ with $y$ 
or $y^{-1}$): the conditions of the Ping-Pong Lemma are now satisfied.

In practice, when attempting to find non-abelian
free subgroups of matrix groups, we use the following variant. 

\begin{lemma} \label{pingpong} 
Let $F$ be a field and let $x,y \in \GL(d,F)$. 
For each formal letter $a \in \{x,x^{-1},y,y^{-1}\}$, suppose there exist 
subsets $\Omega_a$ and $\Omega'_a$ in projective space $P=\P(F^d)$ such 
that the following hold: 
\begin{enumerate} [\rm (i)]
\item $(P\setminus \Omega'_a)a \subset \Omega_a$ for 
each $a\in \{x,x^{-1},y,y^{-1}\}$; 
\item 
$\Omega_b \cap \Omega'_a = \varnothing$ for $a,b \in \{x,x^{-1},y,y^{-1}\}$ 
with $a \neq b^{-1}$;
\item 
$\Omega_a \not \subset \cup_{b\neq a} \Omega_b$.
\end{enumerate}
Then $x$ and $y$ form a free pair.
\end{lemma}

\begin{proof}  Pick $p_a \in \Omega_a \setminus \cup_{b \neq a} \Omega_b$. Suppose $w(x,y)$ is a non-trivial reduced word. By the symmetry of our assumptions, without loss of generality, we may assume that it ends with $x$. Let $a$ be the first letter, so $w=a\cdots x$. If $a \in \{x,x^{-1}\}$, then pick $p=p_y$; 
otherwise, pick $p=p_{x^{-1}}$. In both cases, $p\notin \Omega'_a$. We claim that $p \cdot w(x,y) \in \Omega_x$, implying that $w(x,y)\neq 1$. The assumptions imply that we can follow the trajectory of $p$ letter by letter from one $\Omega_b$ to the next until it lands in $\Omega_x$. But $p \notin \Omega_x$. So $w(x,y) \neq 1$.
\end{proof}

We may apply Lemma \ref{Tits-version}
to $X=\Omega_x \cup \Omega_{x^{-1}}$ and $Y=\Omega_y \cup \Omega_{y^{-1}}$ 
and any $p \notin \bigcup_{a \in \{x,x^{-1},y,y^{-1}\}} (\Omega_a \cup \Omega'_a)$ provided the latter set is not 
all of $P$. 
Let $a\in \{x,x^{-1},y,y^{-1}\}$.
In practice, $\Omega_a$ is a small neighborhood 
$\{p \mid \dist(p,U_a)<\epsilon\}$ around the top eigendirection 
$U_a$, and 
$\Omega'_a$ is a small neighborhood 
$\{p \mid \dist(p,U'_a)\leq \epsilon\}$ around the complementary 
$a$-invariant hyperplane $U'_a$. 
If the regions $\Omega_a,\Omega_a'$ are found 
for each $a$ 
and conditions (i)--(iii) of Lemma \ref{pingpong} hold, then 
$P$ and the $\Omega_a,\Omega_a'$ define a \emph{ping-pong table} for 
$x$ and $y$.

\section{Bi-proximal elements of $\GL(d,\C)$}
\label{sec_biprox}

Given $A \in \GL(d, \C)$, we sort its eigenvalues so that
\[ \Lambda(A) = |\la_1(A)| \geq |\la_2(A)| \geq \ldots \geq |\la_{d-1}(A)| 
    \geq |\la_d(A)|\]
and define the {\em gaps} of $A$ as
\[ \alpha_i(A) = |\la_i(A)|/|\la_{i+1}(A)| \in \R 
   \mbox{ for } 1 \leq i \leq d-1.\]
Let
\[ \alpha(A) = \min \{ \alpha_1(A), \alpha_{d-1}(A) \}.\]
Where it is unambiguous, we write $\la_i = \la_i(A)$,   
$\alpha_i = \alpha_i(A)$, 
and $\alpha = \alpha(A)$.
With this notation, $A$ is {\em proximal} if $\alpha_1(A) > 1$ 
and it is
{\em bi-proximal} if $\alpha(A) > 1$; equivalently, $A$ is bi-proximal if
$A$ and $A^{-1}$ are proximal, since the eigenvalues of $A^{-1}$ are
$\la_d(A)^{-1}, \la_{d-1}(A)^{-1}, \ldots, \la_1(A)^{-1}$.
We record some elementary properties of bi-proximal elements.

\begin{lemma}
\label{biprox}
Let $A \in \GL(d, \C)$ be bi-proximal.
\begin{itemize}
\item[\rm (a)]
If $A \in \SL(d, \C)$, then $|\la_1(A)| > 1$ and $|\la_d(A)| < 1$.
\item[\rm (b)]
$A^{-1}$ is bi-proximal and $\Lambda(A^{-1}) = |\la_d(A)^{-1}|$.
\end{itemize}
\end{lemma}

\subsection{Existence of bi-proximal elements}

To construct a non-abelian free subgroup of a matrix group, 
we must understand when bi-proximal elements can be found.  

\begin{defn} $G\leq \GL(d,\C)$ acts \emph{strongly irreducibly} on 
$\C^d$ if $G$ does not preserve a finite union of non-zero
proper subspaces of $\C^d$.
\end{defn}

This is equivalent to requiring that the Zariski-connected component $G^0$ of 
the identity of the Zariski-closure of $G$ acts irreducibly on $\C^d$. 
If $G \leq \SL(d,\C)$, then this is equivalent to requiring that 
$G^0$ is semisimple and acts irreducibly.
The following is \cite[Proposition 3.11]{tits}.

\begin{prop} 
\label{Sirred}
Suppose $G\leq \GL(d,\C)$ acts strongly irreducibly on $\C^d$.
If $G$ contains a proximal element, then $G$ contains a bi-proximal element; 
indeed, bi-proximal elements form a Zariski-dense subset of $G$.
\end{prop}

The existence of proximal elements of $G$ is a more subtle question.
A compact subgroup of $\GL(d,\C)$ contains no proximal elements, since 
all eigenvalues of its elements have absolute value~$1$. 
Since $\SU(d)$ is compact and Zariski-dense in $\SL(d,\C)$, 
Zariski-density in $\SL(d,\C)$ does not guarantee the existence of 
proximal elements.  However, over $\R$ the following theorem follows from
\cite{abels-margulis-soifer,benoist}.

\begin{theorem}
\label{Zden}
A strongly irreducible subgroup $G$ of $\GL(d, \R)$ contains a proximal 
element if and only if its Zariski-closure  does. 
In particular, if $G$ is Zariski-dense in $\SL(d, \R)$, 
then $G$ contains bi-proximal elements.
\end{theorem}

Given bi-proximal $x, y \in \GL(d,\C)$ 
which satisfy certain distance relations, we show in the remainder
of this section how to 
determine $n \in \N$ so that $x^n$ and $y^n$ form a free pair. 
The main result is Theorem~\ref{Ebound}. 

\subsection{The bound}
A generalized eigenspace 
of $x \in \GL(d, \C)$ is $\ker ( (x-\la_i(x) I_d)^d )$. 
It is an $x$-invariant subspace
on which $x$ acts as the sum of a scalar matrix and a nilpotent matrix.
Let $x$ be bi-proximal, 
so $|\la_1(x)|>|\la_2(x)|$ and $|\la_{d-1}(x)|>|\la_d(x)|$. 

\begin{itemize}
\item[$\bullet$]
Let $U_x$ be the eigenspace corresponding to $\la_1(x)$: 
namely, $U_x = \ker(x-\la_1(x)I_d)$.
\item[$\bullet$]
Let $U_x'$ be the sum of the generalized eigenspaces corresponding
to the eigenvalues $\la_2(x), \ldots, \la_d(x)$: namely, 
$U_x' = \im(x-\la_1(x)I_d)$.
\item[$\bullet$] 
Let $u_x$ be an eigenvector corresponding to $\la_1(x)$; 
since $U_x$ is 1-dimensional, $U_x = \langle u_x \rangle$. 
\end{itemize}

Recall from Lemma \ref{biprox} that $x^{-1}$ is
bi-proximal. Given 
bi-proximal $x, y \in \GL(d,\C)$, let $\Gamma := \{x, \x, y, \y\}$, 
define for $a, b \in \Gamma$ the distance 
$\delta_{a,b} = \dist( u_a, U_b') \in \R$, and set
\begin{equation}\label{delta-def}
 \delta(x,y) := \min \{ \delta_{a,b} \mid a, b \in \Gamma,\ b \neq a^{-1} \}.
\end{equation}
Then $\delta(x, y)$ is the minimum of the $12$ distances from each of the 
four top eigenvectors to each of the three complementary hyperplanes 
other than the one associated to its inverse. 
For $a \in \Gamma$, define 
\[ E_a(x,y) 
 :=  2 \frac{ (d-2) \log\max\{1, (\alpha_1(a)^{1/2} - 1)^{-1}\|a\|_o|
   \la_2(a)|^{-1}\} - 2 \log( \delta(x, y)^2/2) }{\log(\alpha_1(a))}.
\]
Now set 
\begin{equation}\label{E-def}  
E(x,y) := \max \{ E_a(x,y) \mid a \in \Gamma \}.
\end{equation}

The following theorem underpins the bi-proximal method.
It is a quantified version of \cite[Proposition 3.12]{tits}.
The remainder of this section is dedicated to its proof.

\begin{theorem}
\label{Ebound}
Let $x, y \in \GL(d,\C)$ be bi-proximal with $\delta(x,y) > 0$. 
If $n \geq E(x,y)$, then $x^n$ and $y^n$ form a free pair.
\end{theorem}

\subsection{Investigating bi-proximal elements}
\label{lemmas}

As a first step towards the proof of Theorem \ref{Ebound}, we prove an 
optimized version of \cite[Lemma 4.4]{BBr25}, 
which is  based on \cite[Lemma 4.6]{Bre08}. We write $\la_i = \la_i(x)$ and
$\alpha_i = \alpha_i(x)$ throughout.

\begin{lemma}\label{mainlem} 
Let $x \in \GL(d,\C)$  with $\alpha_1 > 1$.  Then  
$\forall u \in \C^d \setminus \{0\}$, $\forall \kappa>0$, $\forall n \ge 1$ \[ \dist(u x^n, u_x)\cdot \dist(u, U_x') \cdot \dist(u_x, U_x')^2
\leq \left(\max\{1,\kappa^{-1}|\la_2|^{-1}\|x\|_o\}\right)^{d-2} 
\left((1+\kappa)\frac{|\la_2|}{|\la_1|} \right)^n.\]
\end{lemma}

Since $\dist$ takes values in $[0,1]$, the left hand side 
is at most $1$. The right hand side tends to $0$ as $n \to \infty$
if and only if $(1+\kappa)|\la_2|/|\la_1| < 1$.
To prove the lemma, we need a preliminary step.

\begin{lemma}\label{conjlem} 
Let $a \in \GL(d,\C)$ and $\kappa>0$. Write $s(a):=\|a\|_o/\Lambda(a)$.
There exists $h \in \SL(d,\C)$ such that $\|h^{-1}a h\|_o\leq (1+\kappa) 
\Lambda(a)$ and 
$$\max\{\|h\|_o,\|h^{-1}\|_o\} 
  \leq \max\left\{1,\frac{s(a)}{\kappa}\right\}^{(d-1)/{2}}.$$ 
\end{lemma}
\begin{proof} 

Conjugating $a$ by a unitary matrix preserves 
$\|a\|_o$ and $\Lambda(a)$, so we may
assume that $a$ is upper triangular.

For $t \in \mathbb{C}^*$, define 
$d_t = \mathrm{diag}(t^{-1}, \ldots, t^{-d}) \in \GL(d,\mathbb{C})$
and set $\tilde{a}(t) := d_t a d_t^{-1}$.
Since $a$ is upper triangular, the $(i,j)$ entry of 
$\tilde{a}(t)$ is $t^{j-i}a_{ij}$,
which vanishes for $j < i$ and is a polynomial in $t$ for $j \geq i$.
Thus $\tilde{a}(t)$ extends to a matrix-valued polynomial in $t$, 
holomorphic on all
of $\mathbb{C}$, where $\tilde{a}(0) = \mathrm{diag}(a_{11}, \ldots, a_{dd})$.
In particular, $\|\tilde{a}(0)\|_o = \max_i |a_{ii}| = \Lambda(a)$.
For $t \neq 0$, define $h_t := t^{(d+1)/2} d_t \in \SL(d,\mathbb{C})$.
Observe that $h_t a h_t^{-1} = d_t a d_t^{-1} = \tilde{a}(t)$ since the 
scalar factors cancel.
Note that, if $|t| \leq 1$, then 
\[
    \max\{\|h_t\|_o, \|h_t^{-1}\|_o\} = \|h_t\|_o = |t|^{-(d-1)/2}.
\]

We make the following claim: for all $|t| \leq 1$,
\begin{equation}\label{eq:claim}
    \|\tilde{a}(t)\|_o \leq \Lambda(a) \cdot s(a) \cdot \frac{1 + s(a)|t|}{s(a) + |t|}.
\end{equation}
The lemma follows immediately by taking $h := h_t^{-1}$, where
$t=\max\{1,{s(a)}/{\kappa}\}^{-1}$.

We now prove the claim.
Fix unit vectors $u, v \in \mathbb{C}^d$ and define $f(t) := \langle u\tilde{a}(t), v \rangle$.
Since $\tilde{a}(t)$ is a polynomial in $t$, the function $f$ is holomorphic on $\mathbb{C}$.
On the unit circle $|t| = 1$, the matrix $d_t$ is unitary, hence
$\|\tilde{a}(t)\|_o = \|d_t a d_t^{-1}\|_o = \|a\|_o$,
so $|f(t)| \leq \|a\|_o$ for all $|t| = 1$.
By the maximum modulus principle, $|f(t)| \leq \|a\|_o$ for all $|t| \leq 1$,
so $g := f/\|a\|_o$ is a holomorphic map from the closed unit disc to the closed
unit disc.
Since $|f(0)| \leq \|\tilde{a}(0)\|_o = \Lambda(a)$, we deduce that 
$|g(0)| \leq \Lambda(a)/\|a\|_o = 1/s(a)$.
Applying the Schwarz--Pick Lemma 
(see, for example, \cite[\S 5.5.2]{Krantz}) to $g$, we deduce that 
\[
    |g(t)| \leq \frac{|g(0)| + |t|}{1 + |t|\,|g(0)|}
    \qquad \text{for all } |t| \leq 1.
\]
The right hand side is increasing in $|g(0)|$, so substituting $|g(0)| \leq 1/s(a)$
and multiplying through by $\|a\|_o$ gives
\[
    |f(t)| \leq \|a\|_o \cdot \frac{1/s(a) + |t|}{1 + |t|/s(a)}
    = \Lambda(a) \cdot s(a) \cdot \frac{1 + s(a)|t|}{s(a) + |t|}.
\]
Taking the supremum over all unit vectors 
$u, v$ establishes~\eqref{eq:claim}.
\end{proof}

\noindent
Recall that $\Lambda(a) \leq ||a||_o$.
Taking $\kappa=1$, this lemma specializes to the following.
\begin{corollary} Let $a \in \GL(d,\C)$. There exists  
$h \in \SL(d,\C)$ such that $\| h^{-1}ah\|_o\leq 2 \Lambda(a)$ and 
$\max\{\|h\|_o,\|h^{-1}\|_o\} \leq s(a)^{(d-1)/{2}}$, 
where $s(a)=\|a\|_o/\Lambda(a)$.
\end{corollary}

\noindent 
We now prove Lemma \ref{mainlem}.

\noindent 
\begin{proof}
 We may decompose $u\in \C^{d}\setminus\{0\}$ as $u=\mu u_x+u'$ 
where $\mu=\mu(u) \in \C$, and 
$u_x$ is an eigenvector of $x$ with eigenvalue $\la_1$ and unit length,
and $u' \in U'_x$. 
We may assume that both $\mu$ and $u'$ are non-zero: 
otherwise, $ux^n \in \C u_x$ and the 
left hand side is zero. The statement follows by multiplying the following 
inequalities
\begin{align}
\dist\left(ux^{n},u_x\right) & \leq\frac{\|u'\|}{\|ux^{n}\|} \max\{1,\kappa^{-1}s(x_{\mid U'_x})\}^{d-2}  ((1+\kappa)|\la_2|)^n  \label{eq:dynamics-1}\\
\dist\left(u_x,U'_x\right) & \leq\frac{\|ux^{n}\|}{\| \mu u_x x^{n}\|}\label{eq:dynamics-2}\\
\dist\left(u_x,U'_x\right) & \leq\frac{\|u\|}{\|u'\|}\label{eq:dynamics-3}\\
\dist\left(u,U'_x\right) & \leq  \frac{\| \mu u_x x^{n}\|}{\|u\|}|\la_1|^{-n}\label{eq:dynamics-4}
\end{align}
and by noting that 
\begin{align}
s(x_{\mid U'_x}) & = \|x_{\mid U'_x}\|_o/\Lambda(x_{\mid U'_x})
                   =\|x_{\mid U'_x}\|_o/|\la_2| \leq  
\|x\|_o |\la_2|^{-1}.
\end{align}

We now prove that the inequalities hold for $n \geq 1$. Observe that 
\begin{align*}
\dist\left(ux^{n},u_x\right) & =\frac{\|ux^{n}\wedge u_x\|}{\|ux^{n}\|}=\frac{\|u'x^{n}\wedge u_x\|}{\|ux^{n}\|}\leq\frac{\|u'x^{n}\|}{\|ux^{n}\|}.
\end{align*}
If $v_x$ is a unit normal vector to $U'_x$,  
then, noting that $u_xx^n \in \C^* u_x$,
\[
\dist\left(ux^{n},U'_x\right)=\frac{|ux^{n}\cdot v_x|}{\|ux^{n}\|}=\frac{|\mu u_x x^{n}\cdot v_x|}{\|ux^{n}\|}= \dist\left(u_x,U'_x\right)\frac{\|\mu u_x x^{n}\|}{\|ux^{n}\|},
\]
so \eqref{eq:dynamics-2} follows. By applying Lemma \ref{conjlem} to 
the restriction of $x$ to $U'_x$, we deduce that 
\begin{align*}
\|u'x^{n}\| &\leq  \|x^{n}_{\mid U'_x}\|_o \|u'\| \\
& \leq \|h^{-1}\|_o \|h\|_o \|h^{-1} x_{\mid U'_x}^n h\|_o\|u'\| \\ &\leq\|h^{-1}\|_o \|h\|_o  \|(h^{-1} x_{\mid U'_x} h)^n\|_o  \|u'\|\\ &\leq \|h^{-1}\|_o \|h\|_o  \|h^{-1} x_{\mid U'_x} h\|_o^n  \|u'\|\\ &\leq \|h^{-1}\|_o \|h\|_o  ((1+\kappa)|\la_2|)^n  \|u'\|\\ &\leq \max\{1,\kappa^{-1}s(x_{\mid U'_x})\}^{d-2}   ((1+\kappa)|\la_2|)^n  \|u'\|,
\end{align*}
so \eqref{eq:dynamics-1} follows. Also
$$\dist(u,U'_x) = \frac{|u \cdot v_x|}{\|u\|} = \frac{|\mu u_x \cdot v_x|}{\|u\|} \leq  \frac{|\mu|}{\|u\|}$$
and
$$\| \mu u_x x^n\|= |\la_1|^n |\mu|,$$
so \eqref{eq:dynamics-4} follows.
Finally, $$\dist(u_x,U'_x)\leq \dist(u_x,u')=\frac{\|u_x\wedge u'\|}{\|u'\|}=\frac{\|u_x \wedge u\|}{\|u'\|}\leq\frac{\|u\|}{\|u'\|},$$
so \eqref{eq:dynamics-3} follows.
\end{proof}

\subsection{The proof of Theorem \ref{Ebound}}

\noindent 
\begin{proof} 
Let $\delta = \delta(x,y)$ and choose $\epsilon = \delta / 2$. This 
guarantees that  the open $\epsilon$-ball around $u_x$ is disjoint 
from the closed $\epsilon$-neighborhood around $U'_x$, $U'_y$ 
and $U'_{y^{-1}}$, and a corresponding disjointness condition holds for 
each of $x^{-1}$, $y$ and $y^{-1}$. We show 
that these eight sets (four open $\epsilon$-balls around 
$u_x,u_{x^{-1}}, u_y$ and $u_{y^{-1}}$ and four closed $\epsilon$-neighborhoods 
of the complementary hyperplanes) define a ping-pong table for 
$x^n$ and $y^n$ for sufficiently large~$n$.

Let $a \in \Gamma$. We write 
$\la_1 = |\la_1(a)|$, $\la_2 = |\la_2(a)|$ and $\alpha =  \alpha_1(a)
= \la_1 / \la_2$.  
Condition (ii) of Lemma~\ref{pingpong} holds by the choice of $\epsilon$.
For condition (iii), note that $u_b \in U'_{b^{-1}}$ for each
$b \in \Gamma$, since $u_b$ is an eigenvector of $b^{-1}$ for its
eigenvalue of smallest absolute value. Hence, for $b \neq a$,
$\dist(u_a, u_b) \ge \dist(u_a, U'_{b^{-1}}) \ge \delta = 2\epsilon$,
so $u_a$ lies in none of the other three $\epsilon$-balls.
By Lemma~\ref{pingpong}, it remains to verify condition (i): namely, 
$a^n$ maps the complement of the closed $\epsilon$-neighborhood around
$U'_a$ into the open $\epsilon$-ball around $u_a$.

For $u \in \C^d \setminus \{0\}$ and
$n \in \N$, assume that $\dist( u, U_a' )>\epsilon$.
The definition of $\delta$ (see \eqref{delta-def}) implies that 
$\dist( u_a, U_a' ) \ge \delta$. If $\dist(ua^n, u_a) = 0$, then
$ua^n$ certainly lies in the open $\epsilon$-ball around $u_a$, so we
may assume that $\dist(ua^n, u_a) > 0$. Applying Lemma \ref{mainlem}, we
deduce that
\begin{align*}
\dist(u a^n , u_a )  \epsilon \delta^2
  &< \dist(u a^n , u_a )  \dist( u, U_a' ) \dist( u_a, U_a' )^2 \\
  &\leq \left(\max\{1,\kappa^{-1}|\la_2|^{-1}\|a\|_o\}\right)^{d-2}
         \left((1 + \kappa)  \alpha_1^{-1} \right)^n \\
  &=: A(a,\kappa, n).
\end{align*}
For $a \in \Gamma$, set 
\[ E_a(\kappa,x,y) := 
    \frac{(d-2) \log( \max \{1, \kappa^{-1} \la_2^{-1} ||a||_o \} )
          - \log( \delta^4 / 4 )}
      {| \log( (1+\kappa) \alpha^{-1}) |}. \]

\begin{lemma}
\label{new}
For each $a \in \Gamma$, choose $\kappa_a>0$ such that
$(1 + \kappa_a)\alpha_1(a)^{-1} < 1$, and assume that
$n \geq E_a(\kappa_a, x, y)$ for each $a \in \Gamma$.
Then $x^n$ and $y^n$ form a free pair.
\end{lemma}

\begin{proof}
Fix $a \in \Gamma$ and write $\kappa = \kappa_a$ and $\alpha = \alpha_1(a)$.
Consider the following equivalences.
\begin{align*}
& A(a,\kappa, n)/(\epsilon \delta^2) \leq \epsilon \\
\iff\quad & A(a, \kappa, n) \leq \epsilon^2 \delta^2 \\
\iff\quad &
 \left(\max\{1,\kappa^{-1}\la_2^{-1}\|a\|_o\} \right)^{d-2} \left((1 + \kappa) \alpha^{-1}\right)^n
  \leq \delta^4 / 4 \\
\iff\quad &
\left((1 + \kappa) \alpha^{-1} \right)^n 
  \leq (\delta^4 / 4) / \left(\max\{1,\kappa^{-1}\la_2^{-1}\|a\|_o\} \right)^{d-2} \\
\iff\quad &
n \log ((1 + \kappa) \alpha^{-1}) \leq \log( \delta^4 / 4) 
       - (d-2) \log (\max\{1,\kappa^{-1}\la_2^{-1}\|a\|_o\}) \\
\iff\quad &
n \geq \frac{(d-2) \log (\max\{1,\kappa^{-1}\la_2^{-1}\|a\|_o\}) 
              - \log( \delta^4 / 4)}
            { | \log ((1 + \kappa) \alpha^{-1}) | }.
\end{align*}
The statement follows.
\end{proof}

Consider $E_a(\kappa,x,y)$. Note that $a \in \Gamma$ and $\delta = \delta(x,y)$ 
depend on $x$ and $y$ only. The parameter $\kappa$ needs only to
satisfy $(1+\kappa)\alpha^{-1} < 1$, where $\alpha = \alpha_1(a)$. 
For fixed $a, x$ and $y$, we view $E_a(\kappa,x,y)$ as a function of
$\kappa$ only. 
For each $a \in \Gamma$, choose $\kappa_a = \alpha^{1/2}-1$. 
Since $\kappa_a$ is fixed, 
we write $E_a(x,y)$ for $E_a(\kappa_a, x,y )$
and conclude that 
\[ E_a(x,y) = 
    2 \frac{(d-2) \log( \max \{1, 
(\alpha^{1/2}-1)^{-1} ||a||_o \la_2^{-1} \} )
          - 2 \log( \delta^2 / 2 )}
      { \log( \alpha ) }. \] 
This completes the proof of Theorem \ref{Ebound}.
\end{proof}

\section{Conjugators in $\GL(d, \C)$}
\label{sec_conjug}

We describe a variation of the method in Section 
\ref{sec_biprox}: instead of searching for two bi-proximal  elements 
of $G \leq \GL(d, \C)$, we search for one bi-proximal 
element $x$ and a suitable {\em conjugator} $b$ for it.  

\begin{defn}\label{defn-conj}
An element $b$ of $G$ is a {\em conjugator} for $x$ if it satisfies 
the following conditions:
\begin{itemize}
\item 
$(u_x)b \not \in U_x' \cup U_{x^{-1}}';$
\item 
$(u_{x^{-1}})b \not \in U_x' \cup U_{x^{-1}}';$
\item 
$u_x \not \in (U_x')b \cup (U_{x^{-1}}')b;$
\item 
$u_{x^{-1}} \not \in (U_x')b \cup (U_{x^{-1}}')b$.
\end{itemize}
\end{defn}

The following is a variation of \cite[Proposition 5]{Dix88}.

\begin{proposition}
\label{dixon}
If $G \leq \GL(d, \C)$ acts strongly irreducibly on $\C^d$ and 
$x \in G$ is bi-proximal, 
then a conjugator $b$ for $x$ exists in $G$.
\end{proposition}

Given a bi-proximal element $x$ and a conjugator $b$ for $x$, we define
$y = x^b = b^{-1}xb$. Note that $x$ and $y$ have the same eigenvalues.
Let $u_y:=(u_x)b$ be an eigenvector of $y$ 
corresponding to the eigenvalue $\la_1$,
and let $U'_y:=(U'_x)b$ be the associated complementary $y$-invariant
hyperplane. Similarly, $u_{y^{-1}}=(u_{x^{-1}})b$ and
$U'_{y^{-1}}=(U'_{x^{-1}})b$. 
Using the function $E_a$ of Section \ref{sec_biprox}, set
\[ B(x,b) := \max\{E_a(x,y) \mid a \in \Gamma \}.\] 
\begin{theorem}\label{Bbound-alternative}
Let $x, b \in G \leq \GL(d,\C)$ where $x$ is bi-proximal and $b$ is a conjugator
for $x$.
If $n \geq B(x,b)$, then $x^n$ and $(x^b)^n$ form a free pair.
\end{theorem}

\begin{proof} 
Let $y = x^b$. By Theorem~\ref{Ebound}, it suffices to show that
$\delta(x,y) > 0$; that is, each of the $12$ distances in
\eqref{delta-def} is positive. Since $\la_1(a)$ is a simple eigenvalue
of $a$, we observe that $u_a \notin U'_a$, so $\delta_{a,a} > 0$ for each
$a \in \{x, \x, y, \y\}$. The four conditions of
Definition~\ref{defn-conj} 
show that the remaining eight distances are positive.
\end{proof}

Theorem~\ref{Bbound-alternative} is a special case of Theorem~\ref{Ebound}.
Nonetheless, the conjugator method is useful, since searching for a bi-proximal
element and a conjugator can yield a pair with a smaller bound than searching 
for two bi-proximal elements independently.

\section{Singular value decomposition in $\GL(d,\R)$ and $\GL(d, \C)$}
\label{sec_cartan}

In this section, we present another method, inspired by 
\cite{breuillard-gelander}, which produces a free pair for 
$G \leq \GL(d, \R)$. The main idea is to consider the singular 
value decomposition of matrices 
and exploit their contracting behaviour on projective space. We describe
the method for matrices over $\R$, but it applies verbatim
to matrices over $\C$, where we use the conjugate transpose $g^*=\ol{g^T}$
of a matrix $g$ instead of its transpose $g^T$.

\subsection{Notation}

If $g \in \GL(d, \R)$, then $g g^T$ is 
positive definite and symmetric and thus has positive real eigenvalues. 
Let $\la_1(g g^T) \geq \ldots \geq \la_d(g g^T)$ denote the eigenvalues
and $\sigma_i(g) = (\la_i(g g^T))^{1/2}$ the singular values. Define
\begin{align}
  \beta(g) &= \frac{\sigma_1(g)}{\sigma_2(g)}, \label{bet0} \\
  \rho(g) &= \min\{\beta(g),\beta(g^{-1})\}. 
\end{align}
Now $\rho(g)$ is the \emph{contraction ratio} of $g$. 
Note that $\sigma_1(g)=||g||_o$, the operator norm of $g$, 
and $\sigma_1(g)\sigma_2(g)$ is the operator norm of the second wedge 
power $\Lambda^2 g$ acting on $\Lambda^2 \R^d$. Define 
\[ \Lip(g)=(\|g\|_o \cdot \|g^{-1}\|_o)^2=(\sigma_1(g)/\sigma_d(g))^2.\]
Using the definition of distance from Section \ref{sec_biprox:distance}, we 
deduce the following.

\begin{lemma}
Every $z \in \GL(d,\R)$ is $\Lip(z)$-Lipschitz on $\mathbb{P}(\R^d)$:
namely, for every $u,v \in \mathbb{P}(\R^d)$,
$$\dist(uz,vz) \leq \Lip(z) \cdot \dist(u,v).$$ 
\end{lemma}

\begin{definition}
The {\em singular value decomposition} of $g \in \GL(d, \R)$ is 
$g=k_1ak_2$ where $k_1$ and $k_2$ are 
orthogonal matrices
 and $a$ is a diagonal 
matrix with diagonal $(\sigma_1(g), \ldots, \sigma_d(g))$.
\end{definition}

\noindent 
Let $v_g$ be an eigenvector of $g^Tg$ corresponding to the largest 
eigenvalue 
and let $v_{g^T}$ be an eigenvector of $gg^T$ corresponding
to the largest eigenvalue. 
If $e_1 = (1, 0, \ldots, 0)$, then 
we may take $v_g = e_1k_2$ and $v_{g^T} = e_1k_1^{-1}$. 
Let $H_g := 
\{ v \in \R^d \mid v \cdot v_{g^T} = 0\}$
be the hyperplane orthogonal to $v_{g^T}$.

\subsection{Choosing generators for a free subgroup}

As a first step, we choose $g \in G$ with large contraction 
ratio $\rho (g)$. Let 
$$\epsilon := \epsilon(g) = \rho(g)^{-1/2}.$$

For $z \in \GL(d, \R)$, 
\begin{align*}
\dist(v_gz, H_g) 
   & = \frac{|v_gz\cdot v_{g^T}|}{\|v_gz\|\cdot\|v_{g^T}\| } \\ 
 \dist(v_{g^{-1}}, H_{g^{-1}}z) 
   & = \frac{|v_{g^{-1}}\cdot v_{(g^{-1})^T}(z^{-1})^T|}{\|v_{{g}^{-1}}\|
       \cdot \|v_{(g^{-1})^T}(z^{-1})^T\|}.
\end{align*}
As a second step, we choose $z \in G$ so that $v_gz$ is as far as 
possible from $H_g$, and $v_{g^{-1}}$ is as far as possible from 
$H_{g^{-1}}z$.  In particular, we require the following:
\begin{align}
\label{fdist0} \dist(v_gz, H_g) &\ge 2\Lip(z)\epsilon \\
\label{sdist0} \dist(v_{g^{-1}}, H_{g^{-1}}z)  &\ge  2\Lip(z)\epsilon.
\end{align}
 
Define $x = gz$.  Let 
$a_x=v_gz$, $B_x=H_g$, $a_{\x}=v_{g^{-1}}$ and $B_{\x}=H_{g^{-1}}z$. 
Now $x$ and $x^{-1}$ have good contracting properties: 
$x$ maps the complement of a small neighborhood around $B_x$ into a 
small ball around $a_x$, and $x^{-1}$ behaves similarly
with respect to $B_{\x}$ and $a_{\x}$. 

It remains to determine $y$.
We obtain $y$ as a conjugate $x^b$ by a suitable conjugator $b$.
Choose $b \in G$ so that $a_xb$ and $a_{\x}b$ are as far as 
possible from $B_x$ and $B_{\x}$. In particular, we require the following:
\begin{align}
\label{aB0} \dist(a_xb,B_xb) &\ge 2\Lip(b)\Lip(z) \epsilon \\
\label{aBbar0} \dist(a_{\x}b,B_{\x}b) &\ge 2\Lip(b)\Lip(z) \epsilon  \\
\label{sepp0} \min_{u \in \{a_x,a_{\x}\}, V \in \{B_x, B_{\x}\}} 
           \{\dist(ub,V), \dist(u,Vb)\}  &\ge 2\Lip(b)\Lip(z) \epsilon.
\end{align}

Note that \eqref{aB0} and \eqref{aBbar0} are stronger than \eqref{fdist0} 
and \eqref{sdist0} respectively, so $\epsilon$ must be small enough to 
accommodate these requirements. 

We now prove the following: if both the contraction ratio $\rho(g)$ and 
the separation between the attracting points and repelling hyperplanes 
are sufficiently large, then $x$ and $y = x^b$ play 
ping-pong.  More precisely, we prove the following theorem.

\begin{theorem}\label{SVDtheorem0}
If $g, z, b \in G$ satisfy $(\ref{aB0})-(\ref{sepp0})$,
then $x = gz$ and $y = x^b$ form a free pair.
\end{theorem}

\subsection{Proof of Theorem \ref{SVDtheorem0}}
The proof proceeds via the following lemmas. Recall that $\epsilon =
\epsilon(g)$ throughout.

\begin{lemma}
\label{contract0} 
Each $g\in \GL(d,\R)$ maps the complement of the $\epsilon$-neighborhood 
of $H_g$ into the $\epsilon$-ball around $v_g$: that is, $\forall u \in 
\mathbb{P}(\R^d),$ $\dist(u,H_g)>\epsilon$ implies $\dist(ug,v_g)<\epsilon$.
\end{lemma}
\begin{proof} 
Assume $g$ has singular value decomposition $k_1 a k_2$ and write $e_1, \ldots,
e_d$ for the standard basis of $\R^d$. 
Let $w=uk_1$ and $w_1 = w\cdot e_1$ and note that $e_1^\perp=\Span(e_2,\ldots,e_d)$.
Now $$\dist(u,H_g)=\dist(u,e_1^\perp k_1^{-1})
=\dist(w,e_1^\perp)=|w_1|/\|w\|>\epsilon.$$ 
Also
$$\dist(ug,v_g)= \dist(uk_1ak_2,e_1k_2)=\dist(wa,e_1)
=\frac{\|wa\wedge e_1\|}{\|wa\|}.$$ 
But $\|wa\|^2=\|wa\wedge e_1\|^2+|\sigma_1w_1|^2$. 
We deduce that 
$$\dist(ug,v_g)^2
= \frac{\|wa\|^2-|\sigma_1w_1|^2}{\|wa\|^2} \leq 
\frac{ \left(|\sigma_2| \cdot \|w\| \right)^2}{|\sigma_1w_1|^2} 
= \left(\frac{|\sigma_2|}{|\sigma_1|}\right)^2 \left(\frac{\|w\|}{|w_1|} \right)^2 
= \beta(g)^{-2} \dist(u,H_g)^{-2}.$$
Hence $\dist(ug,v_g) <  \beta(g)^{-1} \epsilon ^{-1} \leq \epsilon$.
\end{proof}

\begin{lemma}
\label{claim10}
The element
$x=gz$ sends the complement of the closed $\epsilon$-neighborhood of 
$B_x$ into the open $\Lip(z) \epsilon$-ball around $a_x$, which 
itself lies in the intersection of the complement of the closed 
$\epsilon$-neighborhood of $B_x$, the complement of the 
closed $\Lip(b)\epsilon$-neighborhood of $B_xb$, and the complement of 
the closed $\Lip(b)\Lip(z)\epsilon$-neighborhood of $B_{\x}b$.
\end{lemma}
\begin{proof} 
Let $u$ be such that $\dist(u,B_x)>\epsilon$, 
so $$\dist(ux,a_x)=\dist(ugz,v_gz)\leq \Lip(z)\dist(ug,v_g).$$ 
But $\dist(u,B_x)=\dist(u,H_g)>\epsilon$, so $\dist(ug,v_g)<\epsilon$ 
by Lemma \ref{contract0}. Thus $\dist(ux,a_x)< \Lip(z)\epsilon$. 
If $v$ is such that $\dist(v,a_x)<\Lip(z)\epsilon$, then 
$$\dist(v,B_x)\ge \dist(a_x,B_x)-\dist(v,a_x)>\dist(a_x,B_x)-\Lip(z)\epsilon\ge \Lip(z)\epsilon$$ 
provided $\dist(a_x,B_x) \ge 2\Lip(z)\epsilon$, which is the case by \eqref{fdist0}. 
Moreover $$\dist(v,B_xb) \ge \dist(a_x,B_xb)-\dist(v,a_x)>\dist(a_x,B_xb)-\Lip(z)\epsilon \ge \Lip(b)\Lip(z)\epsilon$$ 
provided $\dist(a_x,B_xb) \ge 2\Lip(b)\Lip(z)\epsilon$, which is the case by \eqref{sepp0}. Finally, $$\dist(v,B_{\x}b) \ge \dist(a_x,B_{\x}b)-\dist(v,a_x)>\dist(a_x,B_{\x}b)-\Lip(z)\epsilon \ge \Lip(b)\Lip(z)\epsilon$$ 
provided $\dist(a_x,B_{\x}b) \ge 2\Lip(b)\Lip(z)\epsilon$, which is the case by \eqref{sepp0}.
\end{proof}

\begin{lemma}
\label{claim20}
The element
$x^{-1}=z^{-1}g^{-1}$ sends the complement of the closed 
$\Lip(z)\epsilon$-neighborhood of $B_{\x}$ into the open
$\epsilon$-neighborhood  of $a_{\x}$.
\end{lemma}
\begin{proof} Let $u$ be such that $\dist(u,B_{\x})>\Lip(z)\epsilon$, 
so $\dist(u,H_{g^{-1}}z)>\Lip(z)\epsilon$, 
and thus $\dist(uz^{-1},H_{g^{-1}})>\epsilon$. 
Hence, by Lemma \ref{contract0},  
\begin{equation*}
\dist(ux^{-1},a_{\x})=\dist(uz^{-1}g^{-1},v_{g^{-1}})< \epsilon. \qedhere
\end{equation*}
\let\qed\relax
\end{proof}

\begin{lemma}
\label{claim30}
The element
$y$ sends the complement of the closed $\Lip(b)\epsilon$-neighborhood of 
$B_xb$ into the open $\Lip(b)\Lip(z) \epsilon$-neighborhood  of $a_xb$.
\end{lemma}
\begin{proof} 
Let $u$ be such that $\dist(u,B_xb)>\Lip(b)\epsilon$, so 
$$\dist(uy,a_xb)=\dist(ub^{-1}xb,a_xb)\leq \Lip(b)\dist(ub^{-1}x,a_x).$$ 
But $\dist(u,B_xb)\leq \Lip(b)\dist(ub^{-1},B_x)$, 
so $\dist(ub^{-1},B_x)>\epsilon$. Hence, by Lemma \ref{claim10},
$$\dist(ub^{-1}x,a_x)< \Lip(z)\epsilon,$$ 
so $\dist(uy,a_xb)< \Lip(b)\Lip(z)\epsilon$. 
\end{proof}

\begin{lemma}
\label{claim40}
The element 
$y^{-1}$ sends the complement of the closed 
$\Lip(z)\Lip(b)\epsilon$-neigh\-bor\-hood 
of $B_{\x}b$ into the open $\Lip(b)\epsilon$-neighborhood of $a_{\x}b$.
\end{lemma}
\begin{proof} Let $u$ be such that $\dist(u,B_{\x}b)>\Lip(z)\Lip(b)\epsilon$, 
so 
$$\dist(uy^{-1},a_{\x}b)=\dist(ub^{-1}x^{-1}b,a_{\x}b)\leq \Lip(b)\dist(ub^{-1}x^{-1},a_{\x}).$$ 
But $\dist(u,B_{\x}b)\leq \Lip(b)\dist(ub^{-1},B_{\x})$, 
so $\dist(ub^{-1},B_{\x})>\Lip(z)\epsilon$. Hence, by Lemma \ref{claim20},
$\dist(ub^{-1}x^{-1},a_{\x})< \epsilon$, 
so $\dist(uy^{-1},a_{\x}b)< \Lip(b)\epsilon$. 
\end{proof}

\noindent
We now prove Theorem \ref{SVDtheorem0}.
\begin{proof} 
For $a \in \{x,x^{-1},y,y^{-1}\}$, let $\Omega_a$ be the open ball around 
$a_x,a_{\x},a_xb$ and $a_{\x}b$ of radius  $\Lip(z)\epsilon$, $\epsilon$, $\Lip(b)\Lip(z)\epsilon$, and $\Lip(b)\epsilon$ respectively. Similarly, let 
 $\Omega_a'$ be the closed neighborhoods around $B_x,B_{\x},B_xb$ and $B_{\x}b$ of width $\epsilon$, $\Lip(z)\epsilon$, $\Lip(b)\epsilon$ and $\Lip(b)\Lip(z)\epsilon$ respectively.

When $a=x$, Lemma \ref{claim10} implies that $a$ sends the complement of 
$\Omega_a'$ into $\Omega_a$ and $\Omega_a$ is disjoint from $\Omega'_b$ for $b\neq a^{-1}$. Similarly, Lemmas \ref{claim20}, \ref{claim30} and \ref{claim40} show that $a$ maps the complement of $\Omega_a'$ into $\Omega_a$ for $a \in \{x^{-1},y,y^{-1}\}$. It remains to show that $\Omega_a$ is disjoint from 
$\Omega'_b$ for 
$a \in \{x^{-1},y,y^{-1}\}$
and $b \neq a^{-1}$.
This is done as in the 
case $a=x$ using \eqref{aB0}--\eqref{sepp0}. 
This ensures conditions $(i)$ and $(ii)$ of Lemma \ref{pingpong}. Condition $(iii)$ follows from \eqref{aB0}--\eqref{sepp0}, noting that $a_x \in B_{x^{-1}}$ and $a_{x^{-1}}\in B_x$. We omit the details.

Hence the regions 
$\Omega_a$, $\Omega'_a$ for $a \in \{x,x^{-1},y,y^{-1}\}$ define a
ping-pong table for $x$ and $y$ as in Lemma \ref{pingpong}, 
so $x$ and $y$ form a free pair.
\end{proof}

\section{Groups defined over a non-archimedean local field}
\label{sec_untra}

In this section, we describe how the bound of Section \ref{sec_biprox}
must be modified when $G \leq \GL(d, K)$ for a non-archimedean local 
field $K$ with absolute value $|\cdot |$. Let $\mathcal{O}$ be the ring 
of integers of $K$; that is, $\mathcal{O}=\{x \in K \mid |x|\leq 1\}$. 
Recall that the absolute value $|\cdot|$ on $K$ maps $K^\times$ to 
$\{ \beta^z \mid z \in \Z\} \subset \R_{>0}$ for some $0<\beta<1$.
For example, if $K=\Q_p$ the field of $p$-adic numbers, 
then $\mathcal{O}=\Z_p$ the $p$-adic integers,
$\beta=p^{-1}$, and 
$|p| = p^{-1}$. The valuation $v(x)$ of $x \in K$ satisfies
$|x|=\beta^{v(x)}$. 

The absolute value on $K$ extends in a unique way to an absolute value 
of a finite extension $K'$ of $K$. In fact, it extends uniquely 
to an algebraic closure $\ol{K}$ of $K$. If $e$ is the ramification 
index of 
$K'/K$, then $e \mid [K':K]$ and $v(x) \in 
\frac{1}{e}\Z$, so $v(x) \in \Q$ for $x \in \ol{K}$.

\subsection{Conjugating into $\GL(d, \mathcal{O})$}

To produce a pair of elements of $\GL(d,K)$ that play ping-pong in projective space over $K$, a necessary condition is  that at least some elements of the group $G$ have eigenvalues of modulus different from $1$. 
The following theorem shows that 
to check that this condition holds, it is enough to test it on 
words of length at most $d^2$. 
If the test fails, then $G$ conjugates into $\GL(d,\mathcal{O})$ and thus every element of $G$ has all its eigenvalues of modulus $1$. 

\begin{theorem}\label{bounded}
If $G\leq \GL(d,K)$ is generated by a bounded symmetric set $S$ containing $1$
and each eigenvalue of every $x \in S^{d^2}$ has modulus $1$, 
then $G$ conjugates into $\GL(d, \mathcal{O})$.
\end{theorem}

\begin{proof} This is a simple consequence of 
\cite[Theorem 7]{breuillard-JSR}. 
Equip $K^d$ with the sup-norm $\|v\| = \max_i |v_i|$. Now
$\GL(d,\mathcal{O})$ is exactly the set of matrices $g \in \GL(d,K)$ with
$\|g\|_o = \|g^{-1}\|_o = 1$.
So \cite[Theorem 7]{breuillard-JSR} finds $g \in \GL(d,K')$, where $K'$ is a finite extension of $K$ such that $gGg^{-1} \subset \GL(d,\mathcal{O}')$. Hence $G$ is bounded. But then $L:=\sum_{h \in G} \mathcal{O}^d h$ is a finitely generated $\mathcal{O}$-module of full rank, which is stable under $G$. So $G$ lies in the stabilizer of $L$, hence  it conjugates into $\GL(d,\mathcal{O})$.
\end{proof}

\begin{remark}
\label{choosep-padic}
    If the characteristic of the residue field of $K$ is strictly greater than $d$, then the hypotheses of the previous theorem hold whenever $\Tr(x) \in \mathcal{O}$ for every $x \in S^{d^3}$. This follows immediately from the Newton identities applied to the eigenvalues of each $x\in S^{d^2}$.
\end{remark}

\subsection{Generalized bi-proximal elements}

There is no analogue to Theorem \ref{Zden} over $K$;
see \cite{benoist,quint} for examples of unbounded Zariski-dense groups 
in $\SL(d,\Q_p)$ with no proximal elements. 
Instead, we consider \emph{generalized} bi-proximal elements. 

For $x \in \GL(d,K)$ with eigenvalues $\la_1(x),\ldots,\la_d(x)$ 
ordered so that 
\[ \Lambda(x):= |\la_1(x)| \ge \ldots \ge |\la_d(x)|, \]
let $i_0 \in [1,d/2]$ be 
such that $|\la_{i_0}(x)|>|\la_{i_0+1}(x)|$. 
Define $\alpha_{i_0}(x) = |\la_{i_0}(x)| / |\la_{i_0+1}(x)|$. 
Let $A(x)$ be the sum of the generalized eigenspaces of $x$ 
with eigenvalues $\la_i(x)$ for $i\leq i_0$ and let $R(x)$ be the complementary 
sum of the generalized eigenspaces of $x$ with eigenvalues $\la_i(x)$ for 
$i>i_0$. Where it is unambiguous, 
we write $\la_{i} = \la_{i}(x)$ and $\alpha_i = \alpha_i(x)$.

We consider $x \in \GL(d,K)$ with the property that not only 
$|\la_{i_0}|>|\la_{i_0+1}|$, but the symmetrical condition holds 
for 
$x^{-1}$ (namely, 
$|\la_{d-i_0}| > |\la_{d-i_0+1}|$). If $d$ is even and $i_0=d/2$, then the 
second condition is automatically satisfied when the first is. Such 
elements are {\it generalized bi-proximal}.
Every conjugate of a
generalized bi-proximal element of 
$\GL(d, K)$ is also generalized bi-proximal.

\subsection{Distances and the bound}

Recall from the proof of Theorem \ref{bounded}
that the norm of $v \in K^d$ is the maximal absolute value of its 
coordinates. If $A = (a_{i, j}) \in M(d,\ol{K})$, then 
$||A||_o=\max_{i,j} |a_{i,j}|.$

As in Section \ref{sec_biprox:distance}, 
we define the projective distance between 
$u,v \in K^d \setminus \{0\}$ by 
\begin{equation}\label{dist-def-na} 
\dist(u,v)=\frac{\|u \wedge v\|}{\|u\| \cdot \|v\|}.
\end{equation}
By \cite[2.8.18]{bombieri-gubler}, it satisfies the ultrametric 
triangle inequality: namely,  
\[ \dist(u,v) \leq \max\{\dist(u,w),\dist(w,v)\}. \]

We generalize this definition to subspaces as follows.
Let $U$ and $V$ be subspaces of $\ol{K}^d$. If $u_1,\ldots,u_{k_U}$ 
and $v_{1},\ldots,v_{k_V}$ form bases of $U$ and $V$ respectively, then, 
setting ${\bf u}=u_1 \wedge \cdots \wedge u_{k_U}$ and 
${\bf v}=v_{1} \wedge \cdots \wedge v_{k_V}$,  we define
\begin{equation*} 
\dist(U,V) = \frac{\|\bf{u} \wedge \bf{v}\|}{\|\bf{u}\| \cdot \|\bf{v}\|}. 
\end{equation*}
Note that $\dist(U,V)>0$ if and only if $U \cap V = \{0\}$. It can be 
shown (see \cite[\S 3]{BBr25}) that for $u \in U \setminus \{0\}$ 
\begin{equation} 
\label{incl} 
\dist(u,V) := \inf_{v \in V \setminus \{0\}} \dist(u,v) \ge \dist(U,V).
\end{equation}
It follows from this and the ultrametric triangle inequality that if 
$0<\epsilon\leq \dist(U,V)$, then the open $\epsilon$-neighborhoods of $U$ 
and $V$ are disjoint. 
Indeed, if $\dist(w,U)<\epsilon$ and $\dist(w,V)<\epsilon$, 
then there exist $u\in U$ and $v\in V$ such that $\dist(w,u)<\epsilon$ and 
$\dist(w,v)< \epsilon$ and hence $\dist(u,v)<\epsilon$. 
We conclude 
that $\dist(U,V)<\epsilon$ by \eqref{incl}.

Let $x \in \GL(d, K)$ be generalized 
bi-proximal and let
$y = x^b$ be a conjugate of $x$.
Note that $A(x)b = A(y)$ and $R(x)b = R(y)$ and also $A(\xx)b = A(\yy)$ 
and $R(\xx)b = R(\yy)$.
Consider the following distances: 
\begin{itemize}
\item 
$\dist( A(x), R(x))$ and $\dist(A(\xx), R(\xx))$;
\item 
$\dist( A(x)b, R(x)b)$ and $\dist(A(\xx)b, R({\xx})b)$;
\item 
$\dist( A(x)b, R(x))$ and $\dist(A(x), R(x)b)$;
\item 
$\dist( A(x)b, R(\xx))$ and $\dist(A(x), R({\xx})b)$;
\item 
$\dist( A(\xx)b, R(x))$ and $\dist(A(\xx), R(x)b)$;
\item 
$\dist( A(\xx)b, R(\xx))$ and $\dist(A(\xx), R({\xx})b)$.
\end{itemize}
Set $\delta(x,x^b)$ to be the minimum of these 12 distances. 

\begin{theorem}
\label{padic}
Let $K$ be a non-archimedean local field and let $G$ be an unbounded Zariski-dense subgroup of $\SL(d,K)$. Then $G$ contains a generalized bi-proximal element. Moreover, if $x$ is such an element, then there exists $b \in G$ with $\delta(x, x^b) > 0$. 
\end{theorem}
\begin{proof} The unboundedness of $G$ and Theorem \ref{bounded} (applied to a boundedly generated unbounded subgroup of $G$) imply that at least some $g \in G$ has not all its eigenvalues of modulus $1$. Since $g \in \SL(d,K)$, this implies that there is $i_0$ with $1\leq i_0<d$ such that $|\lambda_{i_0}(g)|>|\lambda_{i_0+1}(g)|$. In other words, $G$ contains a proximal element for the  fundamental irreducible 
representation $\Lambda^{i_0}(K)$. The arguments of 
\cite[Proposition 3.11]{tits} show how to
produce a bi-proximal element from a proximal one in 
an irreducible representation. This proves the existence of $x$. 
If subspaces $U$ and $V$ of $\ol{K}^d$ have complementary 
dimension, 
then the set of $g \in \SL(d,\ol{K})$ such that 
$\dist(Ug,V)>0$ is a \emph{non-empty} Zariski-open set. 
Indeed, $\SL(d,\ol{K})$ acts transitively on subspaces of dimension $k$ 
for each $k$.  Therefore each 
finite intersection of such sets is again non-empty and Zariski-open, hence it
must contain a Zariski-dense set of elements from $G$. This guarantees the existence of $b \in G$ with $\delta(x,x^b)>0$.
\end{proof}

We say that $b$ 
is a {\it conjugator} for $x$ if $\delta(x,x^b)>0$.
Let $\delta_a:=\dist(A(a),R(a))$ and 
$\ol{N}(a):=(\|a\|_o \cdot \|a^{-1}\|_o)^{d-2}$.
Define 
\[ \ell(a, b):= \frac{\log( \ol{N}(a)/(\delta (x, x^b) \cdot \delta_a)^2)}{\log(\alpha_{i_0}(a))}.\]
Fix $x$, $b$ and $y := x^{b}$ as in Theorem \ref{padic}.
Now set 
\[ H_{\mathrm {opt}}(x,b):=\max\{ \ell(a, b) \mid a \in \{x, x^{-1}, y, y^{-1}\} \}.\]
In the remainder of this section we prove the following theorem.
\begin{theorem}\label{Bboundultra}
Let $K$ be a non-archimedean local field. Let $x, b \in G\leq \GL(d,K)$ 
where $x$ is generalized bi-proximal and $b$ is a conjugator
for $x$. 
If $n > H_{\mathrm {opt}}(x,b)$, 
then $x^n$ and $(x^b)^n$ form a free pair. 
\end{theorem}

\begin{remark}\label{rem-Nbar}
Note that $\ol{N}(x)=\ol{N}(x^{-1})$ and $\ol{N}(y)=\ol{N}(y^{-1})$
and $\alpha_{i_0}(x^{-1})=\alpha_{d-i_0}(x)$. 
We may get a better bound if we replace $\ol{N}(a)$ by 
$N(i_0,a)$, defined in Lemma \ref{mainlemultra} below, for each 
$a \in \{x, x^{-1}, y, y^{-1}\}$. The theorem remains 
valid with this replacement.
\end{remark}

For $x \in \GL(d,K)$ let $v(x)$ be the 
valuation of $x$: namely, the integer such that 
$\beta^{v(x)}=\max_{i,j} |x_{i,j}|$. 
Note that $v(x)=\min_{i,j} v(x_{i,j})$.
Abusing notation, we also denote by $v(\delta)$ the rational number 
such that $\delta=\beta^{v(\delta)}$. 
Each of $v(x)$,  $v(x^{-1})$, $v(y)$ and $v(y^{-1})$ is integral. 
Each subspace $A(a)$ and $R(a)$ used to define $\delta$ is the sum of
the generalized eigenspaces for all eigenvalues in a range of absolute
values; such a sum is defined over $K$, and hence $v(\delta)$ is integral too.

The following corollary is immediate. It gives a 
slightly worse bound in a more compact form. 

\begin{corollary} With $x,b$ as above, $x^n$ and $(x^b)^n$ form
a free pair for each $n > H(x,b)$ where 
\[H(x,b) = 
\frac{(d-2) \left(\max\{|v(x)|+|v(x^{-1})|,|v(x^b)|+|v((x^b)^{-1})|\} \right)
+4|v(\delta)|}
     {\min\{|v(\la_{i_0})-v(\la_{i_0+1})|,|v(\la_{d-i_0})-v(\la_{d-i_0+1})|\} }. \]
\end{corollary}

\subsection{Investigating generalized bi-proximal elements}

We state and prove analogues to the results in Section 
\ref{lemmas}.

\begin{lemma}
\label{mainlemultra} 
Let $x \in \GL(d,\ol{K})$  with 
$\alpha_{i_0} = |\la_{i_0}| / |\la_{i_0+1}|>1$. 
Then  $\forall u \in \ol{K}^d \setminus \{0\}$, $\forall n \ge 1$ \[ \dist(u x^n, A(x))\cdot \dist(u, R(x)) \cdot \dist(A(x), R(x))^2 
\leq N(i_0,x)  \cdot \alpha_{i_0}^{-n}\]
where 
\[ N(i_0,x)= \left(\frac{\|x_{\mid R(x)}\|_o}{|\la_{i_0+1}|}\right)^{d-i_0-1}\cdot  \left(\frac{\|x_{\mid A(x)}^{-1}\|_o}{|\la_{i_0}|^{-1}}\right)^{i_0-1}.\]
\end{lemma}
\noindent
Note that $$N(i_0,x) \leq (\|x\|_o \cdot \|x^{-1}\|_o)^{d-2}.$$
To prove the lemma, we need a preliminary step. 
\begin{lemma}
\label{conjlemultra} 
Let $a \in \GL(d,\ol{K})$. There exists $h \in \SL(d,\ol{K})$ such that $\|h^{-1}a h\|_o\leq  \Lambda(a)$ and 
$$\max\{\|h\|_o,\|h^{-1}\|_o\} \leq  s(a)^{\frac{d-1}{2}},$$ 
where $s(a):=\|a\|/\Lambda(a)$.
\end{lemma}
\begin{proof} The proof is the same as for Lemma \ref{conjlem} {\it mutatis 
mutandis}, taking $t\in \ol{K}$ such that $|t|=s(a)^{-1}$. We can assume that
$a$ is upper triangular, then we conjugate by the diagonal matrix
$d_t = \diag(t^{-1}, \ldots, t^{-d})$ as in that proof. We omit the details. 
\end{proof}

\noindent
We now prove Lemma \ref{mainlemultra}.

\noindent
\begin{proof}
We may decompose $u\in \ol{K}^{d}\setminus\left\{ 0\right\} $ as $u=u_{A}+u_{R}$ where $u_A \in A(x)$ and $u_R \in R(x)$. 
The statement follows by multiplying the following inequalities: 
\begin{align}
\dist\left(ux^{n},A(x)\right) & \leq\frac{\|u_{R}\|}{\|ux^{n}\|}\left(\frac{\|x_{\mid R(x)}\|}{\Lambda(x_{\mid R(x)})}\right)^{d-i_0-1} |\la_{i_0+1}|^{n}\label{eq:dynamic-1}\\
\dist\left(A(x),R(x)\right) & \leq\frac{\|ux^{n}\|}{\|u_{A}x^{n}\|}\label{eq:dynamic-2}\\
\dist\left(A(x),R(x)\right) & \leq\frac{\|u\|}{\|u_{R}\|}\label{eq:dynamic-3}\\
\dist\left(u,R(x)\right) & \leq\frac{\|u_{A}x^{n}\|}{\|u\|} \left(\frac{\|x_{\mid A(x)}^{-1}\|}{\Lambda\left(x_{\mid A(x)}^{-1}\right)}\right)^{i_0-1}|\la_{i_0}|^{-n}.\label{eq:dynamic-4}
\end{align}
We now prove that the inequalities hold for $n \geq 1$.
Let ${\bf a}$ (resp.\ ${\bf r}$) be a unit pure wedge vector representing
$A(x)$ (resp.\ $R(x)$). Then
\begin{align*}
\dist\left(ux^{n},A(x)\right) & =\frac{\|ux^{n}\wedge{\bf a}\|}{\|ux^{n}\|}=\frac{\|u_{R}x^{n}\wedge{\bf a}\|}{\|ux^{n}\|}\leq\frac{\|u_{R}x^{n}\|}{\|ux^{n}\|}
\end{align*}
and, using \eqref{incl} in the last inequality,
\[
\dist\left(ux^{n},R(x)\right)=\frac{\|ux^{n}\wedge{\bf r}\|}{\|ux^{n}\|}=\frac{\|u_{A}x^{n}\wedge{\bf r}\|}{\|ux^{n}\|}\geq \dist\left(A(x),R(x)\right)\frac{\|u_{A}x^{n}\|}{\|ux^{n}\|},
\]
so \eqref{eq:dynamic-2} follows. By applying Lemma \ref{conjlemultra} 
to the restriction of $x$ on $R(x)$ and noting
that $\Lambda\left(x_{\mid R(x)}\right)=|\la_{i_0+1}|$, we deduce that 
\[
\|u_{R}x^{n}\|\leq \|u_{R}h\|\cdot \|\left(h^{-1}x h\right)^{n}\| \cdot \|h^{-1}\| \leq \left(\frac{\|x_{\mid R(x)}\|}{\Lambda(x_{\mid R(x)})}\right)^{d-i_0-1}|\la_{i_0+1}|^{n}\|u_{R}\|,
\]
so \eqref{eq:dynamic-1} follows. By applying the same lemma to
the restriction of 
$x^{-1}$ on $A(x)$ and noting that 
$\Lambda\left(x_{\mid A(x)}^{-1}\right)=|\la_{i_0}|^{-1}$, we deduce that 
\[
\|u_{A}\|\leq \left(\frac{\|x_{\mid A(x)}^{-1}\|}{\Lambda\left(x_{\mid A(x)}^{-1}\right)}\right)^{i_0-1}|\la_{i_0}|^{-n}\|u_{A}x^{n}\|\,\,,
\]
so \eqref{eq:dynamic-4} follows since $\dist\left(u,R(x)\right)\leq\frac{\|u_{A}\|}{\|u\|}$.
Finally, by \eqref{incl}, 
$$\dist\left(A(x),R(x)\right)\leq \dist\left(A(x),u_{R}\right)=\frac{\|{\bf a}\wedge u_{R}\|}{\|u_{R}\|}=\frac{\|{\bf a}\wedge u\|}{\|u_{R}\|}\leq\frac{\|u\|}{\|u_{R}\|},$$
so \eqref{eq:dynamic-3} follows.
\end{proof}

\subsection{The proof of Theorem \ref{Bboundultra}}

\noindent
\begin{proof} We apply Lemma \ref{mainlemultra} to each of $x$, 
$x^{-1}$, $y=x^b$ and $y^{-1}$. 
Recall that $\delta_x=\dist(A(x),R(x))$ 
and $\delta$ is the minimum of the 
12 distances listed before Theorem \ref{padic}.
By the ultrametric triangle inequality,  
the open $\delta$-neighborhood of $A(x)$ is disjoint from 
the open $\delta$-neighborhoods of $R(x)$, $R(y)$ and $R(y^{-1})$ and 
a corresponding disjointness condition holds for each 
of $x^{-1}$, $y$ and $y^{-1}$.  

By Lemma \ref{mainlemultra}, 
if $u \in \ol{K}^d \setminus \{0\}$ satisfies $\dist(u,R(x)) \geq   \delta$, 
then 
\[ \dist(u x^n, A(x)) \leq  \frac{N(i_0,x) \cdot \alpha_{i_0}^{-n}}
        {\delta \cdot \delta_x^2}.\]  
The latter is less than $\delta$ if 
\[ \alpha_{i_0}^{-n}
    < \frac{(\delta \cdot \delta_x)^2}{N(i_0,x)},
 \] 
or equivalently if
\begin{equation}\label{bestvalue} n > \ell'(x):= \frac{\log( N(i_0,x)/(\delta \cdot \delta_x)^2 )}{\log(\alpha_{i_0})}.\end{equation}
If 
\[ n > \max\{ \ell'(x), \ell'(x^{-1}), \ell'(y), \ell'(y^{-1})\},\]
then 
Lemma \ref{mainlemultra} implies that $x^n$ and $y^n$ play ping-pong 
on projective space $\mathbb{P}(\ol{K}^d)$ with attracting neighborhoods 
the open $\delta$-neighborhoods $A(x)_\delta$, $A(x^{-1})_\delta$, $A(y^{-1})_\delta$ and $A(y)_\delta$ around $A(x)$, $A(x^{-1})$, $A(y^{-1})$ and $A(y)$ respectively, and repelling neighborhoods the open $\delta$-neighborhoods 
around $R(x)$, $R(x^{-1})$, $R(y^{-1})$ and $R(y)$. 
Hence $x^n$ sends the complement of  $R(x)_\delta$ (and, in particular, 
each of $A(x)_\delta$, $A(y^{-1})_\delta$ and $A(y)_\delta$)  
into  $A(x)_\delta$. 
The corresponding conditions (i) and (ii) in Lemma \ref{pingpong} hold 
for $x^{-1}$, $y$ and $y^{-1}$ and we check condition (iii) similarly.
Therefore $x^n$ and $y^n$ form a free pair.
Since $N(i_0,a) \le \ol{N}(a)$, this maximum is at most 
$H_{\mathrm{opt}}(x,b)$. 
\end{proof}

\section{Cartan decomposition over a non-archimedean local field}
\label{sec-padicSVD}

In this section, we present a non-archimedean analogue to the method 
of Section \ref{sec_cartan}. It produces a free pair for 
$G \leq \GL(d, K)$, where $K$ is a non-archimedean local field. 
The main idea is to consider the  Cartan decomposition
$$\GL(d,K)=\GL(d,\mathcal{O})A \GL(d,\mathcal{O}),$$ 
where $A$ is the subgroup of diagonal matrices of the 
form $\diag(\pi^{n_1},\cdots,\pi^{n_d})$ for $n_i\in \Z$ 
and $\pi$ is a uniformizer of $K$. 

\subsection{Notation}

The {\it Cartan decomposition} of $g \in \GL(d, K)$ is 
$g=k_1 a k_2$, 
where $a=\diag(\pi^{n_1},\cdots,\pi^{n_d})$ and 
$k_1,k_2 \in \GL(d,\mathcal{O})$. 
This is the non-archimedean analogue 
of the singular value decomposition.
We may assume that the diagonal entries are ordered so that 
$n_1\leq \ldots \leq n_d$.
This decomposition is the Smith normal form of any 
scalar multiple 
of $g$ with entries in $\mathcal{O}$.  
Let $\sigma_i(g) = |\pi^{n_i}|$ and define
\begin{align}
  \beta(g) &= \frac{\sigma_1(g)}{\sigma_2(g)}, \label{bet} \\
  \rho(g) &= \min\{\beta(g),\beta(g^{-1})\}. \label{rho}
\end{align}
As before, $\rho(g)$ is the contraction ratio of $g$. 
Note that $\sigma_1(g)=||g||_o$, the operator norm of $g$, 
and $\sigma_1(g)\sigma_2(g)$ is the operator norm of the second wedge 
power $\Lambda^2 g$ acting on $\Lambda^2 K^d$. Define 
\[ \Lip(g):=(\|g\|_o \cdot \|g^{-1}\|_o)^2.\]
The relevant distance $\dist$ is defined in \eqref{dist-def-na}
and satisfies the ultrametric triangle inequality 
$\dist(u,v)\leq \max\{\dist(u,w), \dist(w,v)\}$. 
Also $u\cdot v = uv^T= \sum_i u_iv_i \in K.$
\begin{lemma}
Every $z \in \GL(d,K)$ is $\Lip(z)$-Lipschitz on $\mathbb{P}(K^d)$:
namely, for every $u,v \in \mathbb{P}(K^d)$,
$$\dist(uz,vz) \leq \Lip(z) \cdot \dist(u,v).$$ 
\end{lemma}

\noindent 
Let $v_g := e_1k_2$ and $v'_g := e_1k_1^T$, and let
$H_g := \{v \in K^d \mid v\cdot v'_g = 0\}$. 
If $J$ is the involution sending $e_i$ to $e_{d+1-i}$ for
$i \in \{1,\ldots,d\}$, then the Cartan decomposition of $g^{-1}$ 
is $(k_2^{-1} J)(J a^{-1} J)(J k_1^{-1})$.

\subsection{Choosing generators for a free subgroup}

As a first step, we choose $g \in G$ with large contraction 
ratio $\rho (g)$. Let 
$$\epsilon := \epsilon(g) = \rho(g)^{-1/2}.$$
For $z \in \GL(d, K)$, 
\begin{align*}
\dist(v_gz, H_g) 
   & = \frac{|v_gz\cdot v'_{g}|}{\|v_gz\|\cdot\|v'_{g}\| } \\ 
\dist(v_{g^{-1}}, H_{g^{-1}}z) 
   & = \frac{|v_{g^{-1}}\cdot v'_{g^{-1}}z^{-T}|}{\|v_{{g}^{-1}}\|
       \cdot \|v'_{g^{-1}}z^{-T}\|}.
\end{align*}

As a second step, we choose $z \in G$ so that $v_gz$ is as far as 
possible from $H_g$, and 
$v_{g^{-1}}$ is as far as possible from $H_{g^{-1}}z$.
In particular, we require the following:
\begin{align}
\label{fdist} \dist(v_gz, H_g) &\ge \epsilon \\
\label{sdist} \dist(v_{g^{-1}}, H_{g^{-1}}z)  &\ge  \Lip(z)\epsilon.
\end{align}
 
Define $x = gz$.  Let 
$a_x=v_gz$, $B_x=H_g$, $a_{\x}=v_{g^{-1}}$ and $B_{\x}=H_{g^{-1}}z$. 
Now $x$ and $x^{-1}$ have good contracting properties: 
$x$ maps the complement of a small neighborhood around $B_x$ into a 
small ball around $a_x$, and $x^{-1}$ behaves similarly
with respect to $B_{\x}$ and $a_{\x}$. 

It remains to determine $y$.
We obtain $y$ as a conjugate $x^b$ by a suitable conjugator $b$.
Choose $b \in G$ so that $a_xb$ and $a_{\x}b$ are as far as 
possible from $B_x$ and $B_{\x}$. In particular, we require the following:
\begin{align}
\label{aB} \dist(a_xb,B_xb) &> \Lip(b)\Lip(z) \epsilon \\
\label{aBbar} \dist(a_{\x}b,B_{\x}b) &> \Lip(b)\Lip(z) \epsilon  \\
\label{sepp} \min_{u \in \{a_x,a_{\x}\}, V \in \{B_x, B_{\x}\}} 
           \{\dist(ub,V), \dist(u,Vb)\}  &> \Lip(b)\Lip(z) \epsilon.
\end{align}
Note that \eqref{aB} and \eqref{aBbar} are stronger than \eqref{fdist} 
and \eqref{sdist}, so $\epsilon$ must be small enough to 
accommodate these requirements. 

We now prove the following: if both the contraction ratio $\rho(g)$ and 
the separation between the attracting points and repelling hyperplanes 
are sufficiently large, then $x$ and $y =  x^b$ 
form a free pair.
More precisely, we prove the following theorem.

\begin{theorem}\label{SVDtheorem}
If $g, z, b \in G$ satisfy $(\ref{aB}) - (\ref{sepp})$,
then $x = gz$ and $y = x^b$ form a free pair.
\end{theorem}

\subsection{Proof of Theorem \ref{SVDtheorem}}

The proof proceeds via the following lemmas. Recall that $\epsilon =
\epsilon(g)$ throughout.

\begin{lemma}
\label{contract} 
Each $g\in \GL(d,K)$ maps the complement of the closed $\epsilon$-neighborhood 
of $H_g$ into the open $\epsilon$-ball around $v_g$: that is, $\forall u \in 
\mathbb{P}(K^d),$ $\dist(u,H_g)>\epsilon$ implies $\dist(ug,v_g)<\epsilon$.
\end{lemma}
\begin{proof} 
Assume $g$ has Cartan decomposition $k_1 a k_2$.
Let $w=uk_1$ and  $w_1 = w\cdot e_1$ and note that 
$e_1^\perp=\Span(e_2,\ldots,e_d)$. 
Now $$\dist(u,H_g)=|u\cdot v'_g|/ (\|u\| \cdot \|v'_g\|)
=\dist(u,e_1^\perp k_1^{-1})=\dist(w,e_1^\perp)=|w_1|/\|w\|>\epsilon.$$ 
Also
$$\dist(ug,v_g)= \dist(uk_1ak_2,e_1k_2)=\dist(wa,e_1)
=\frac{\|wa\wedge e_1\|}{\|wa\|}.$$ 
Recall that $\sigma_1=|\pi^{n_1}|$, so 
$\|wa\|=\max_i |w_i \pi^{n_i}| \ge |w_1|\sigma_1$. We deduce that  
$$\dist(ug,v_g)
= \frac{\max_{i>1}|\pi^{n_i}w_i|}{\|wa\|}\leq 
\frac{\sigma_2}{\sigma_1} \frac{\|w\|}{|w_1|} \leq  \beta(g)^{-1} \dist(u,H_g)^{-1}.$$
Hence $\dist(ug,v_g) <  \beta(g)^{-1} \epsilon ^{-1} \leq \epsilon$.
\end{proof}

\begin{lemma}
\label{claim1}
The element
$x=gz$ sends the complement of the closed $\epsilon$-neighborhood of 
$B_x=H_g$ into the open $\Lip(z) \epsilon$-neighborhood  of $a_x=v_gz$.
\end{lemma}
\begin{proof} 
Let $u$ be such that $\dist(u,B_x)>\epsilon$, 
so $$\dist(ux,a_x)=\dist(ugz,v_gz)\leq \Lip(z)\dist(ug,v_g).$$ 
But $\dist(u,B_x)=\dist(u,H_g)>\epsilon$, so $\dist(ug,v_g)<\epsilon$ 
by Lemma \ref{contract}; thus $\dist(ux,a_x)< \Lip(z)\epsilon$. 
\end{proof}

\begin{lemma}
\label{claim2}
The element
$x^{-1}=z^{-1}g^{-1}$ sends the complement of the closed
$\Lip(z)\epsilon$-neighborhood of $B_{\x}$ into the open
$\epsilon$-neighborhood  of $a_{\x}$.
\end{lemma}
\begin{proof} Let $u$ be such that $\dist(u,B_{\x})>\Lip(z)\epsilon$, 
so $\dist(u,H_{g^{-1}}z)>\Lip(z)\epsilon$, 
and thus $\dist(uz^{-1},H_{g^{-1}})>\epsilon$. 
Hence, by Lemma \ref{contract}, 
\begin{equation*}
    \dist(ux^{-1},a_{\x})=\dist(uz^{-1}g^{-1},v_{g^{-1}})< \epsilon. \qedhere
\end{equation*}
\let\qed\relax
\end{proof}

\begin{lemma}
\label{claim3}
The element
$y$ sends the complement of the closed $\Lip(b)\epsilon$-neighborhood of 
$B_xb$ into the open $\Lip(b)\Lip(z) \epsilon$-neighborhood  of $a_xb$.
\end{lemma}
\begin{proof} 
Let $u$ be such that $\dist(u,B_xb)>\Lip(b)\epsilon$, so 
$$\dist(uy,a_xb)=\dist(ub^{-1}xb,a_xb)\leq \Lip(b)\dist(ub^{-1}x,a_x).$$ 
But $\dist(u,B_xb)\leq \Lip(b)\dist(ub^{-1},B_x)$, 
so $\dist(ub^{-1},B_x)>\epsilon$. Hence, by Lemma \ref{claim1}, 
\mbox{$\dist(ub^{-1}x,a_x)< \Lip(z)\epsilon$}, 
so $\dist(uy,a_xb)< \Lip(b)\Lip(z)\epsilon$. 
\end{proof}

\begin{lemma}
\label{claim4}
The element $y^{-1}$ sends the complement of the 
closed $\Lip(z)\Lip(b)\epsilon$-neigh\-bor\-hood 
of $B_{\x}b$ into the open $\Lip(b)\epsilon$-neighborhood of $a_{\x}b$. 
\end{lemma}
\begin{proof} Let $u$ be such that $\dist(u,B_{\x}b)>\Lip(z)\Lip(b)\epsilon$, 
so 
$$\dist(uy^{-1},a_{\x}b)=\dist(ub^{-1}x^{-1}b,a_{\x}b)\leq \Lip(b)\dist(ub^{-1}x^{-1},a_{\x}).$$ 
But $\dist(u,B_{\x}b)\leq \Lip(b)\dist(ub^{-1},B_{\x})$, 
so $\dist(ub^{-1},B_{\x})>\Lip(z)\epsilon$. Hence, by Lemma \ref{claim2},
$\dist(ub^{-1}x^{-1},a_{\x})< \epsilon$, 
so $\dist(uy^{-1},a_{\x}b)< \Lip(b)\epsilon$. 
\end{proof}

\noindent 
We now prove Theorem \ref{SVDtheorem}.

\noindent 
\begin{proof}
For $a \in \{x,x^{-1},y,y^{-1}\}$, let $\Omega_a$ be the open ball 
around $a_x$, $a_{\x}$, $a_xb$ and $a_{\x}b$ of radius $\Lip(z) \epsilon$, 
$\epsilon$, $\Lip(b)\Lip(z)\epsilon$ and $\Lip(b)\epsilon$ respectively. 
Similarly, 
let $\Omega_a'$ be the closed neighborhood around 
$B_x$, $B_{\x}$, $B_xb$ and $B_{\x}b$ of width $\epsilon$, 
$\Lip(z)\epsilon$, $\Lip(b)\epsilon$ and $\Lip(b)\Lip(z)\epsilon$ 
respectively. The ten conditions of \eqref{aB}--\eqref{sepp} 
together with the ultrametric triangle inequality for $\dist$ guarantee 
that the disjointness condition $(ii)$ of Lemma \ref{pingpong} holds. 
The last four lemmas show that condition $(i)$ of 
Lemma \ref{pingpong} also holds. Condition $(iii)$ follows from \eqref{aB}--\eqref{sepp}, noting that $a_x \in B_{x^{-1}}$ and $a_{x^{-1}} \in B_x$. 
Hence these regions in $\P(K^d)$ define a ping-pong table, so $x$ and $y$ 
form a free pair. 
\end{proof}

\section{Algorithms and their implementations}
\label{algo_var1}

The theory presented in Sections \ref{sec_biprox}--\ref{sec-padicSVD}
underpins methods to construct a free pair for $G = \langle S \rangle 
\leq \GL(d, F)$ where $F$ is $\C$ or a non-archimedean local field. 
The methods rely on the existence of certain elements of $G$.
We now consider their termination.

\begin{theorem}
\label{termC}
Let $G \leq \GL(d,\C)$.
\begin{enumerate}
\item[\rm (1)]
Suppose that $G$ acts strongly irreducibly on $\C^d$ and contains
a proximal element.
\begin{enumerate}
\item[\rm (a)]
$G$ contains
a bi-proximal element $x$ and a conjugator $b$ with $\delta(x,x^b) > 0$.
\item[\rm (b)]
$G$ contains
bi-proximal elements $x$ and $y$ with $\delta(x,y) > 0$.
\item[\rm (c)]
The singular value
decomposition method terminates.
\end{enumerate}
\item[\rm (2)]
If $G$ is Zariski-dense in $\SL(d,\R)$, then each of the bi-proximal, 
conjugator and singular value decomposition methods terminates. 
\end{enumerate}
\end{theorem}

\begin{proof}
\mbox{}
\begin{enumerate}
\item[\rm (1a)]
Proposition \ref{Sirred} implies that $G$ has a bi-proximal element $x$. 
(In fact, $G$ contains infinitely many bi-proximal 
elements.)  Proposition \ref{dixon} asserts that $G$ contains a conjugator 
$b$ for $x$. Since $x$ and $y = x^b$ are bi-proximal, 
$\delta_{a,a} > 0$ for each $a \in \{x, x^{-1}, y, y^{-1}\}$.
Since $b$ is a conjugator for $x$, the remaining eight relevant
distances are also positive. Hence $\delta(x,y) > 0$.
\item[\rm (1b)]
This follows from (1a): a conjugate of a bi-proximal element is 
bi-proximal.
\item[\rm (1c)]
As in (1a), $G$ contains a bi-proximal element $g$. 
By the spectral radius formula, applied to $g$ and to $\Lambda^2 g$, 
$\beta(g^r)^{1/r} \to \alpha_1(g)$ and 
$\beta(g^{-r})^{1/r} \to \alpha_1(g^{-1}) = \alpha_{d-1}(g)$ 
as $r \to \infty$, so $\rho(g^r)^{1/r} \to \alpha(g) > 1$. 
Hence $\rho(g^r)$ grows exponentially in $r$ or, equivalently, 
$\epsilon(g^r)$ tends to $0$. Thus it suffices 
to find $z$ so that the distances in (\ref{fdist0}) and 
(\ref{sdist0}) are non-trivial and similarly for the conjugator $b$ and
(\ref{aB0})--(\ref{sepp0}), which is guaranteed by (1a).
\item[(2)]
This follows from (1) and Theorem \ref{Zden}. \qedhere
\end{enumerate}
\end{proof}
Detinko, 
Flannery \& Hulpke \cite{DFH-Zariski} describe a procedure that may 
conclude, in certain situations, that $G \leq \SL(d, \Z)$ is 
Zariski-dense in $\SL(d,\R)$.

\begin{theorem}
\label{termK}
Let $G \leq \SL(d, K)$ for a non-archimedean local field $K$ and assume that
$G$ is Zariski-dense and unbounded. 
\begin{enumerate}
\item[\rm (a)]
The group $G$ contains a generalized
bi-proximal element $x$ and a conjugator $b$ with $\delta(x,x^b) > 0$.
\item[\rm (b)]
The non-archimedean conjugator method terminates.
\item[\rm (c)]
If $G$ contains a  bi-proximal element, then  the 
Cartan decomposition method terminates. 
\end{enumerate}
\end{theorem}

\begin{proof}
Part (a) follows directly from Theorem \ref{padic}. Part (b) follows from (a)
and Theorem \ref{Bboundultra}. The proof of (c) is 
similar to that of Theorem \ref{termC}(1c) using (\ref{fdist})--(\ref{sepp}).
\end{proof}

Assume $G$ acts reducibly on $F^d$.
If $W > U$ are $G$-submodules, then there is an induced action of
$G$ on the smaller space $W/U$. Associated with this is a 
group homomorphism
$G \ra H$ with $H \leq \GL(W/U)$. If a free pair
is found for $H$, then arbitrary preimages
form a free pair for $G$. 

Relying on Theorems \ref{termC} and \ref{termK}, and with clear abuse of 
language, we refer to our methods as algorithms. We now discuss their 
implementation in Oscar. 
In the rest of this section, we assume that the inputs to our algorithms
are matrices defined over $\Q$ or an algebraic number field, and so allow
precise computations; we
also use other arithmetic models and present these
in Section \ref{arith} as a preliminary step.  
Our implementations of the archimedean algorithms 
require the number field to have a real embedding, 
and work in a fixed such embedding; the
others accept input over $\Q$.
In Section \ref{implbiprox}, we discuss in detail 
our implementation of the bi-proximal algorithm. 
We comment briefly on the others in Section \ref{implrest}. 

\subsection{Types of arithmetic}
\label{arith}

Our algorithms apply to finite extensions of a number field $F$. 
These can be described by algebraic extensions which allow precise 
arithmetic, but this model is often prohibitively slow. In the following 
we consider alternatives.

\subsubsection{Algebraic arithmetic}
Let $f \in F[t]$ be an irreducible polynomial. Then $F' = F[t]/(f)$ is 
a finite extension of $F$ with precise arithmetic. Note that 
$F' \cong F(\alpha)$
for each root $\alpha$ of $f$ and $[F':F] = \deg(f)$.

\subsubsection{Floating-point arithmetic}
If $F \subseteq \C$, then elements of $F$ are 
$a+ib$ with $a,b \in \R$ and each $w \in \R$ can 
be described by a floating-point approximation $\ol{w}$ so that 
$|w-\ol{w}| < \epsilon$ for an initial precision $\epsilon$. 
Arithmetic for floating-point numbers accumulates rounding errors,
so it is difficult to control the final precision of a 
calculation using this model.

\subsubsection{Interval arithmetic}
A real number $w$ can also be described by a rational interval 
$[m-r, m+r]$ containing $w$. The midpoint $m$ and the radius $r$
of the interval are in $\Q$.
In \cite{AH1983, Johansson2013arb} an arithmetic for intervals is
introduced. This includes elementary arithmetic $(+,-, \times ,/)$ and 
the application of functions such as $\log$. 

\subsubsection{$p$-adic arithmetic}
The arithmetic for elements of extensions of  $\Q_p$ 
is similar to interval arithmetic in the archimedean case: 
the intervals are replaced by closed balls around rational numbers 
with respect to the $p$-adic topology.
See \cite{CG2000, Caruso2017, Doris2021} for details.

\subsection{The bi-proximal algorithm}
\label{implbiprox}

We first discuss how to evaluate the relevant bound $E(x,y)$
for a pair $x, y$ of bi-proximal elements of $\GL(d, \C)$,
and then present a procedure to select $x$ and $y$. 
To verify that a matrix is bi-proximal, we approximate the roots of its
characteristic polynomial (see Section~\ref{evaluate-bound}), increasing the
precision until $|\la_1| > |\la_2|$ and $|\la_{d-1}| > |\la_d|$ are
certified.

\subsubsection{Evaluating the bound}\label{evaluate-bound}

Assume that $x$ and $y$ are bi-proximal elements of $\GL(d, F)$ where $F$ is 
an algebraic number field. Let $\Gamma = \{x, x^{-1}, y, y^{-1}\}$. Recall that
$E(x,y) = \max \{ E_a(x,y) \mid a \in \Gamma \}.$

We split our computation into three procedures. A critical feature of each 
is the arithmetic used: algebraic (and thus precise) or interval arithmetic
(and thus fast). Floating-point arithmetic is not used, as the accumulation
of rounding errors may lead to incorrect results. 

\medskip
\noindent
{\sc Procedure Preprocess $(a \in \Gamma)$}. 
\begin{enumerate}
\item[(1)]
Construct $\Fa = F(\la_1(a))$. Use algebraic arithmetic in $\Fa$ in
the following.
\item[(2)]
Compute an eigenvector $u_a$ corresponding to the eigenvalue $\la_1(a)$;
note that $u_a$ is a non-zero element of $\ker( a - \la_1(a) I_d)$.
\item[(3)]
Compute the sum $U_a'$ of the generalized eigenspaces corresponding to 
the eigenvalues 
$\la_2(a), \ldots, \la_d(a)$; note that $U_a' = \im( a - \la_1(a) I_d)$.
\item[(4)]
Calculate a non-zero vector $c_a \in \ker(U_a'^T)$, so that
$\sum_i (c_a)_i w_i = 0$ for $w \in U_a'$.
\item[(5)]
Return $u_a$ and $c_a$ as vectors in $\Fa^d$.
\end{enumerate}

The field $\Fa = F(\la_1(a))$ is constructed as follows.
Let $\mu_a(t) \in F[t]$ be the minimal polynomial of $a$. We factorize
this over $F$, so $\mu_a(t) = g_1^{e_1}(t) \cdots g_m^{e_m}(t)$ 
for distinct irreducible factors $g_1, \ldots, g_m$. We approximate the roots 
of each $g_i$ (see \cite{Pan1997}) to any desired precision and so 
identify the root $\la_1(a)$ with the largest absolute value. (Similarly, the
root $\la_2(a)$ is determined for later applications.) If $g_i(\la_1(a)) = 0$,
then $\Fa = F[t] / (g_i(t))$.

We use algebraic arithmetic in {\sc Preprocess}, since applications 
of linear algebra to compute kernels and images of matrices are 
rarely robust if done with other arithmetic models. 
Since 
$\Fa$ is an extension of degree at
most $d$, this computation is typically practical.
Since algebraic arithmetic in the next two procedures 
would require computation in an algebraic extension 
of degree up to $d^4$ over $F$, 
we employ interval arithmetic.

\medskip
\noindent
{\sc Procedure ComputeDelta $(x,y)$.} 
\begin{enumerate}
\item[(1)]
Compute $u_a$ and $c_a$ for $a \in \Gamma$ with {\sc Preprocess}, and
translate them to intervals using the embedding of $\Fa$ into $\C$ that
maps $\la_1(a)$ to the eigenvalue of $a$ of largest modulus.
The complex conjugate $n_a$ of $c_a$ is a normal vector to $U_a'$ in 
the sense of Section~\ref{sec_prelim}.
\item[(2)]
Using interval arithmetic,
evaluate $\delta_{a,b} = \dist(u_a, U_b') 
= |u_a \cdot n_b|/ (\|u_a\| \cdot \|n_b\|)$ 
for the $12$ relevant pairs $(a,b)$ with $a,b \in \Gamma$.
\item[(3)]
Using interval arithmetic, compute the minimum 
$\delta(x,y)$ of the values of $\delta_{a,b}$.
\item[(4)] 
If the interval for $\delta(x, y)$ contains $0$, then return false; 
otherwise, return $\delta(x, y)$.
\end{enumerate}

\noindent
If {\sc ComputeDelta} returns false, then we choose a different pair 
$x,y$ and restart. 
Otherwise, we use the following procedure to evaluate $E(x, y)$, 
where we use the observation that there are functions $L(a)$ and $M(a)$ such
that $$E_a(x,y) = L(a) - M(a) \log(\delta(x,y)).$$ 

\medskip
\noindent
{\sc Procedure EBound$(x, y)$.} 
\begin{enumerate}
\item[(1)]
Compute $\delta = \delta(x,y)$ with {\sc ComputeDelta}. If $\delta$ is 
false, then return false.
\item[(2)]
For each $a \in \Gamma$ use interval arithmetic for each of the following.
\begin{itemize}
\item
Compute $|\la_1(a)|$, $|\la_2(a)|$ and $||a||_o$. Recall that $\alpha_1(a)
= |\la_1(a)|/|\la_2(a)|$. 
\item
Evaluate $L(a)$ and $M(a)$ as intervals in $\R$. 
\item
Determine $E_a(x,y) = L(a) - M(a) \log(\delta(x,y))$ as the interval $[m-r,m+r]$.
\item
Let $c(a)$ be the ceiling $m+r$ of the determined interval.
\end{itemize}
\item[(3)]
Return the minimal $n \in \N$ with $n \geq c(a)$ for each $a \in \Gamma$.
\end{enumerate}

\subsubsection{Selecting bi-proximal elements}

In Figure~\ref{selection}, we present the procedure {\sc SelectPair},
which searches for bi-proximal elements $x, y$ of $G = \langle S \rangle$
with a small value of {\sc EBound}$(x,y)$. It takes as input $S$ and three
positive integers: a bound $k$ on word length, a stability parameter $m$,
and a bound $N$ on the number of random selections.

At each step, the procedure constructs two random words $w, v$ of length at
most $k$ in the free group of rank $|S|$, and evaluates them on $S$ to obtain
$x = w(S)$ and $y = v(S)$. If $x$ or $y$ is not bi-proximal, or if
{\sc EBound}$(x,y)$ returns false, then the pair is discarded. Otherwise,
$c :=$ {\sc EBound}$(x,y)$ is compared with the smallest value found so
far. If this minimum has not decreased for
$m$ successive admissible pairs, then the procedure
returns a pair $x, y$ achieving the
minimum, the words $w, v$, and the value of {\sc EBound}$(x,y)$.
We make at most $N$ selections;
if no admissible pair is found, then the procedure
returns  \textsc{fail};
otherwise, if stability is not achieved, then 
the procedure returns the best admissible pair constructed
and its associated bound.

\begin{figure}[ht]
\begin{algorithmic}
  \State $stable \leftarrow 0$;\ $min_b \leftarrow \infty$;\ $good \leftarrow$ \textsc{fail}
  \For{$i = 1$ \textbf{to} $N$}
    \State choose random words $w, v$ of length at most $k$
    \State $(x,y) \leftarrow (w(S), v(S))$
    \If{$x$ or $y$ is not bi-proximal} \textbf{continue} \EndIf
    \State $c \leftarrow$ \textsc{EBound}$(x,y)$
    \If{$c$ is \textsc{false}} \textbf{continue} \EndIf
    \If{$c < min_b$}
      \State $min_b \leftarrow c$;\ $good \leftarrow (x,y,w,v)$;\ $stable \leftarrow 0$
    \Else
      \State $stable \leftarrow stable + 1$
      \If{$stable \geq m$} \State \Return $good, min_b$ \EndIf
    \EndIf
  \EndFor
  \State \Return $good, min_b$ 
\end{algorithmic}
\caption{Procedure \textsc{SelectPair}$(S,k,m,N)$}\label{selection}
\end{figure}

By default, we choose a small value for $k$. If {\sc SelectPair} fails,
then we increase $k$ and call it again. In practice, if $G$ contains
bi-proximal elements, then they are typically easy to find; $d \log d$
is often a useful upper bound for $k$. If $G$ satisfies the hypotheses of
Theorem~\ref{termC}(1), then bi-proximal pairs with $\delta(x,y) > 0$ exist.

\subsection{The other algorithms}
\label{implrest}

Our implementations of 
the algorithms of Sections \ref{sec_conjug}--\ref{sec-padicSVD} use
variations of the machinery described above.

\subsubsection{The conjugator algorithm}
We use an easy modification of {\sc SelectPair} to 
select $x$ and $b$ and then evaluate the bound $E(x, x^b)$ as before.

\subsubsection{Singular value decomposition algorithm}\label{svd-implement}
The input is $G = \langle S \rangle \leq \GL(d, \C)$.
We first search randomly for an element
$g$ with ``large'' contraction ratio $\rho(g)$. As shown in the proof 
of Theorem~\ref{termC}(1c), if $G$ contains a bi-proximal element $g$, 
then $\rho(g^r)$ grows exponentially in $r$, so elements with 
arbitrarily large contraction ratio are readily found.
Let $\epsilon(g) = \rho(g)^{-1/2}$.
Next, we sample for elements $z$ and $b$ of $G$ satisfying
\eqref{aB0}--\eqref{sepp0}. If $\epsilon(g)$ is
small enough, then these are typically easy to find:
usually, the generators of $G$ suffice.

\subsubsection{The non-archimedean conjugator algorithm}\label{padic-implement}
We assume as input $G = \langle S \rangle \leq \GL(d, \Q)$ (to allow
precise computations with $G$). 

We search for a prime $p$ as follows. 
Let $P$ be the set of primes dividing the denominators of entries of matrices
in $S\cup S^{-1}$. If $p\notin P$, then 
$G\leq\mathrm{GL}(d,\mathbb Z_p)$. 
We first test  a small sample of 
``short" random words and select a prime $p\in P$ for 
which a trace is not $p$-integral.
If no prime is found, then we consider all freely reduced words of
length at most $d^2$ in $S\cup S^{-1}$. For each resulting matrix
and each $p\in P$, we test whether its characteristic polynomial
has $p$-integral coefficients and unit constant term, or equivalently
whether all of its eigenvalues have $p$-adic absolute value $1$. A failed
test determines $p$. If all tests pass, then Theorem~\ref{bounded} proves
boundedness at every finite prime, and this method does not apply.

Keeping the selected prime fixed, we now search randomly for a
generalized bi-proximal element $x$ and a conjugator $b$, and apply
Theorem \ref{Bboundultra}.

Eigenvalue valuations are determined exactly from the Newton polygon.
To compute the attracting and repelling subspaces for $x$ and $x^{-1}$,
we multiply the relevant slope factors by their full multiplicities,
evaluate these polynomials at $x$, and take kernels over $\mathbb Q_p$.
These local computations use finite $p$-adic precision, increased
as needed up to a prescribed limit.

\subsubsection{Cartan decomposition algorithm}
The input is $G = \langle S \rangle \leq \GL(d, \Q)$. For each sampled
$g$, we consider the primes $p$ dividing the numerator or denominator of
the trace of $g$; for each such $p$ with $\epsilon(g) < 1$, we 
embed $\Q$ in $\Q_p$ and use the strategy of Section~\ref{svd-implement}.
We compute in $\GL(d, \Q_p)$ the Cartan decomposition of $g \in G$.
By computing the Smith normal form of an integral multiple of $g$, the 
decomposition can be computed exactly.

\section{Verification using cones in $\SL(d, \Q)$}
\label{section_cones}

Assume that $G $ is a finitely generated subgroup of
$\SL(d, \Q)$.  Each of our algorithms returns a pair 
$x, y$ of elements of $G$ and a positive integer $n$ such
that $x^n$ and $y^n$ form a free pair.
We describe a procedure that takes as input 
$x$, $y$, and a positive integer $l$ and, using cones, 
tries to construct a ping-pong table for $x^l$ and $y^l$. 
If it succeeds, then $x^l$  and $y^l$ form a free pair for $G$. 
The procedure typically verifies freeness for
a value of $l$ that is significantly smaller than $n$.
We have implemented the procedure in both {\sc Magma} \cite{magma}
and Oscar. Since it relies only on rational approximations, 
the associated arithmetic is exact.

A {\it cone} is a subset of $\Q^d$ which is the convex hull of a sequence of
rays defined by \emph{generating vectors} $v_1, \ldots, v_k \in \Q^d$.
It equals 
$$\left\{ \sum_{i=1}^k \beta_i v_i : \beta_i \in \Q,\,\beta_i \ge 0 \right\}.$$

The procedure is designed for the case where each of $x$, $x^{-1}$,
$y$ and $y^{-1}$ has a unique eigenvalue $\lambda$ of maximal absolute value,
where $\lambda \in \R$ and $|\lambda| > 1$, and the (real) eigenspace
corresponding to $\lambda$ has dimension $1$. We denote these eigenspaces
by $\langle u_g \rangle$, with $g \in \{x,{x^{-1}},y,{y^{-1}} \}$;
they should be mutually distinct (and preferably not too close together).
Our implementation checks only that $x$ and $y$ each have eigenvalues
of absolute value greater than $1$ and less than $1$. Since all
subsequent tests are exact, the correctness of a positive answer does
not depend on these assumptions.

For each such $g$, we define the cone $C^+_g$ as a rational approximation
of a cone around the axis defined by $u_g$.
More precisely, we first define the rational approximation
$u_g' := v\, g^{50}$ of $u_g$, where $v$ is a short vector with
entries in $\{\pm 1, \pm 2\}$; we found this power-iteration approach
to be the most effective. We then complete
to form a basis $u_g' = b_1,b_2,\ldots,b_d$ of $\Q^d$, 
where the vectors $b_i$ are mutually orthogonal and of length both 
approximately $1$ and at least $1$.  We choose a small positive rational
number $\tau$ and define $C^+_g$ and $C^-_g$ to be the cones with generator 
vectors $\{ b_1 \pm \tau b_i : 2 \le i \le d \}$ and 
$\{ -v : v \in  C^+_g \}$ respectively, 
and set $C^{\pm}_g := C^+_g \cup C^-_g$.

We use the four sets $C^{\pm}_g$ with $g \in \{x,{x^{-1}},y,{y^{-1}} \}$
as the regions for application of the Ping-Pong Lemma. 
We require $l > 0$ such that, for each such $g$, the element $g^l$ maps
each of the sets $C^{\pm}_h$ other than $C^{\pm}_{g^{-1}}$ into $C^{\pm}_g$.
(For example, $y^l$ must map each of $C^{\pm}_x$,
$C^{\pm}_{x^{-1}}$, and $C^{\pm}_y$ into $C^{\pm}_y$.)
To ensure the validity of the application, we also check that neither
$u'_x$ nor $u'_{x^{-1}}$ lies in $C^{\pm}_y \cup C^{\pm}_{y^{-1}}$, and
neither $u'_y$ nor $u'_{y^{-1}}$ lies in $C^{\pm}_x \cup C^{\pm}_{x^{-1}}$.
(If that fails, then we can try reducing $\tau$.)
These conditions suffice. For example, let $w$ be a non-trivial reduced
word ending in $x$. If its first letter is $x^{\pm 1}$ or $y$, then
$u'_y \, w(x^l,y^l) \in C^{\pm}_x$; if its first letter is $y^{-1}$,
then $u'_{y^{-1}} \, w(x^l,y^l) \in C^{\pm}_x$. Neither starting point
lies in $C^{\pm}_x$, so $w(x^l,y^l) \neq 1$. The other cases are similar.
The attracting property of eigenvectors corresponding to eigenvalues with
largest absolute values should ensure that this method works for suitably
small $\tau$, large $l$, and sufficiently accurate rational approximations.

Since a cone is the convex hull of its generating vectors and
the map defined by $g$ is linear, it suffices to check, for each
$h \neq g^{-1}$, that $g^l$ maps either all of the generating vectors
of $C^+_h$ into $C^+_g$, or all of them into $C^-_g$; the
corresponding statement for $C^-_h = -C^+_h$ then follows.
(Since $\tau$ is small, the generating vectors are close together, so we  
expect one of these two possibilities to hold.)

The parameter $\tau$ can be chosen by a user. Its default value is
$10^{-2}$ when $d \le 4$; otherwise, it is $10^{-4}$.
It is not unusual for the verification to fail with the 
default value and succeed with a smaller value.
Elements produced by the SVD and Cartan algorithms are usually very 
strongly contracting, but their attracting directions may lie close to 
the repelling hyperplanes of the other elements. Verifying such a pair 
with $l = 1$ may then require much smaller values of $\tau$; in our 
experiments, values as small as $10^{-18}$ were needed.

The procedure requires membership tests for vectors in cones. 
Initially, we used the {\sc Magma} functionality for cones to do this.
Since it 
involves exact computations with rational numbers, a successful calculation of
this type is theoretically guaranteed to be correct. However, for
larger values of $d$ (depending also on the size of the entries of the
generators of $G$), these membership tests are prohibitively slow.
Our implementation instead tests a sufficient condition for membership,
described below. This is significantly faster, and a positive answer is
guaranteed to be correct; a negative answer may mean only that the
test was inconclusive.

Suppose that we want to test whether $v \in \Q^d$ lies in $C^+_g$. 
Observe that  
each of the generating vectors $b_1 + \tau b_i$ of $C^+_g$ lies in a 
hyperplane $H$ through the vector $b_1$, where $H$ is orthogonal 
to $b_1$. A straightforward calculation, using the fact that each 
vector $b_i$ has length at least $1$, shows that every $u \in H$ 
satisfying $\|u - b_1\| < \tau /\sqrt{d-1}$ is in $C^+_g$.
We check that this condition holds for $u$ equal to the (unique) scalar
multiple of the vector $v$ that 
lies in $H$. To avoid calculations with real
numbers, we check whether $\sum_{i=1}^d \gamma_i^2 < \tau^2/(d-1)$,
where $u-b_1 = (\gamma_1, \ldots,\gamma_d) \in \Q^d$.

\section{The 4-strand Burau representation}
\label{sec_burau}

Consider the following elements of $\SL(3,\Z[t, t^{-1}])$:
\[ A = \begin{pmatrix}
    0  & 0 & -t^{-1} \\ 0 & -t & -t^{-1} + t \\ -1 & 0 & -t^{-1} + 1
\end{pmatrix} \quad\quad B = \begin{pmatrix}
    -t & t & 0 \\ 0 & 1 & 0 \\ 0 & t^{-1} & -t^{-1} 
\end{pmatrix} \]
The question of whether $\langle A, B \rangle$ is a non-abelian free group is 
intimately connected to the 
longstanding problem to decide whether the Burau representation
of the 4-strand braid group is faithful. In fact, the freeness of 
$\langle A,B \rangle$ is equivalent to the
faithfulness of the Burau representation; see~\cite[Theorem 3.19]{Bir1974}. 
Recently Bharathram, Birman \& Brendle \cite{BBB} 
proved that this representation is faithful.

In~\cite{WZ2015} it was shown that $\langle A^3, B^3 \rangle$ is a 
non-abelian free group; see also~\cite{BT2018}. 
Our results give an independent, computer-assisted proof of that
result.
Beridze, Bigelow \& Traczyk \cite{BBT20} also prove that
$\langle A^3,B^3 \rangle$ is a non-abelian free subgroup of 
\mbox{$\SL(3, (\Z/m\Z)[t, t^{-1}])$} for every integer $m >1$.

Note that $A^b = B$ where
\[ b = \begin{pmatrix} 0 & 0 & -1 \\ 1 & 0 & -1 \\ 0 & 1 & -1\end{pmatrix}.\] 

We report on two experiments using 
our implementations \cite{GITPAGE} in Oscar version 1.8.2 on 
a machine with an Intel Ultra 7 155U processor.
By specializing to specific $t = t_0 \in \Q^\times$, we obtain matrices 
$A(t_0), B(t_0) \in \SL(3, \Q)$. We apply the algorithms to the group
$\langle A(t_0), B(t_0) \rangle$.
Each value of $t_0$ that we use is prime, and for the $p$-adic
algorithms we take $p = t_0$.
In Table \ref{tab:burau-new}, for each algorithm we record:
\begin{itemize}
\item the value of $t_0$; 
\item the maximal length $k$ of the words in
$S := \{A(t_0), B(t_0)\}$ representing $x$ and $y$; for all but the
bi-proximal algorithm, $y = x^b$, and we use the length of the freely
reduced word for $b^{-1}xb$;
\item the exponent $\ell \in \N$ certified by the corresponding theorem,
so that $x^{\ell}$ and $y^{\ell}$ form a free pair (for the SVD and Cartan
algorithms, $\ell = 1$);
\item the time in seconds to construct $x, y$  and $\ell$.
\end{itemize}

\begin{table}[h]
  \centering
  \caption{Results for the five algorithms}
  \label{tab:burau-new}
  \renewcommand{\arraystretch}{1.2}
  \begin{tabular}{r *{5}{r r r}}
    \toprule
    & \multicolumn{3}{c}{bi-proximal}
    & \multicolumn{3}{c}{conjugator}
    & \multicolumn{3}{c}{$p$-adic conjugator}
    & \multicolumn{3}{c}{SVD}
    & \multicolumn{3}{c}{{Cartan}} \\
    \cmidrule(lr){2-4}
    \cmidrule(lr){5-7}
    \cmidrule(lr){8-10}
    \cmidrule(lr){11-13}
    \cmidrule(lr){14-16}
    $t_0$
      & Time & $k$ & $\ell$
      & Time & $k$ & $\ell$
      & Time & $k$ & $\ell$
      & Time & $k$ & $\ell$
      & Time & $k$ & $\ell$ \\

    17   & 0.03\,s & 2 & 14 & 0.46\,s & 6 & 15 & 0.10\,s & 7 & 5 & 0.12\,s & 96 & 1 & 1.73\,s & 32 & 1 \\
    101  & 0.02\,s & 2 & 17 & 0.56\,s & 7 & 22 & 0.18\,s & 8 & 8 & 0.29\,s & 88 & 1 & 2.29\,s & 29 & 1 \\
    1009 & 0.03\,s & 2 & 11 & 0.62\,s & 7 & 12 & 0.27\,s & 7 & 8 & 0.54\,s & 94 & 1 & 2.37\,s & 34 & 1 \\

    \bottomrule
  \end{tabular}
\end{table}

The searches are randomised: our implementations use a variant of
{\sc SelectPair} that
minimises the product of the certified bound and the maximal word length,
rather than the bound alone, so $k$ and $\ell$ vary between runs.
Table~\ref{tab:burau-new} records, for each entry, the rounded median
values over $10$ runs. Pairs returned by the SVD and Cartan algorithms were also
verified to be free using the cones method of
Section~\ref{section_cones} with $l = 1$; this required values of
$\tau$ between $10^{-2}$ and $10^{-18}$.

\begin{table}[hb]
  \centering
  \caption{Evaluating bounds and verification}
  \label{tab:burau1-new}
  \renewcommand{\arraystretch}{1.2}
  \begin{tabular}{r *{4}{r r}}
    \toprule
    & \multicolumn{2}{c}{bi-proximal}
    & \multicolumn{2}{c}{conjugator}
    & \multicolumn{2}{c}{$p$-adic conjugator}
    & \multicolumn{2}{c}{cones} \\
    \cmidrule(lr){2-3}
    \cmidrule(lr){4-5}
    \cmidrule(lr){6-7}
    \cmidrule(lr){8-9}
    $t_0$
      & Time & $\ell$
      & Time & $\ell$
      & Time & $\ell$
      & Time & $\ell$ \\
    \midrule
    17 & 1.1\,ms & 14 & 1.1\,ms & 14 & 1.3\,ms & 5 & 84.9\,ms & 3 \\
    101 & 1.2\,ms & 12 & 1.2\,ms & 12 & 1.6\,ms & 5 & 147.2\,ms & 3 \\
    1009 & 1.2\,ms & 11 & 1.2\,ms & 11 & 1.5\,ms & 5 & 268.9\,ms & 3 \\
    \bottomrule
  \end{tabular}
\end{table}

In Table~\ref{tab:burau1-new}, we consider the pair $A(t_0)$ and $B(t_0)$: 
we record the time in milliseconds taken to evaluate 
the bi-proximal bound of Theorem \ref{Ebound}, 
the conjugator bound  of Theorem \ref{Bbound-alternative},  
and the $p$-adic conjugator bound  of Theorem \ref{Bboundultra}. 
We also record the time taken by the cones method of Section 
\ref{section_cones} to verify that $A(t_0)^{\ell}$ and $B(t_0)^{\ell}$ 
form a free pair.

The bounds of Theorems~\ref{Ebound} and~\ref{Bbound-alternative} give 
$\ell = 14, 12, 11$, and that of Theorem~\ref{Bboundultra} gives $\ell = 5$, 
while the cones method verifies directly that $A(t_0)^3$ and $B(t_0)^3$ 
form a free pair. 

\section{Comments on performance} \label{commentary}

We applied the algorithms to random groups
$\langle A, B\rangle \leq \SL(n,\Q)$ for $n \leq 9$, where the entries of
$A$ and $B$ are rationals with bounded numerators and denominators. 
Scripts that reproduce Tables~\ref{tab:burau-new} and~\ref{tab:burau1-new}
and generate such random examples, with user-chosen dimension, entry bound
and number of groups, are included in our GitHub package~\cite{GITPAGE}.  
Our experiments showed that the $p$-adic conjugator algorithm
was the fastest in every dimension, and the Cartan algorithm returned
$\ell = 1$ with short words; the certified values of $\ell$ returned by the
bi-proximal, conjugator and $p$-adic conjugator algorithms grew 
with $n$. The words constructed by the singular value decomposition algorithm
were the longest.

The cones method of Section~\ref{section_cones} can often confirm that
the pairs returned by the bi-proximal and conjugator algorithms 
play ping-pong for a power $l$ much smaller than the certified $\ell$. 
Its cost grows
with the dimension and, more markedly, with the size of the entries of $x$
and $y$, and a suitable value of $\tau$ must often be found by trial and error.
For pairs returned by the singular value decomposition algorithm, whose entries
are typically very large, it succeeded only in dimensions up to about 8,
depending on the size of the entries in the input group.

\bibliographystyle{abbrv}
\bibliography{refs}

\end{document}